\documentclass[11pt]{amsart}

\usepackage{amssymb,amsmath,epsf,bbm,a4wide}

   \newcommand{\ZZ}{\mathbb{Z}}
   \newcommand{\NN}{\mathbb{N}}
   \newcommand{\QQ}{\mathbb{Q}}
   \newcommand{\RR}{\mathbb{R}}
   
   \newcommand{\EE}{\mathbb{E}}
   \newcommand{\II}{\mathbb{I}}
   \newcommand{\gG}{\varGamma}
   \newcommand{\gL}{\varLambda}
   \newcommand{\gS}{\varSigma}

 \newcommand{\den}{\mathrm{den}}

 \newcommand{\bs}[1]{\boldsymbol{ #1}}

 \newtheorem{lemma}{Lemma}
 \newtheorem{prop}{Proposition}
 \newtheorem{theorem}{Theorem}
 \newtheorem{coro}{Corollary}
 
 \theoremstyle{definition}
  \newtheorem{defin}{Definition}

\numberwithin{equation}{section}
\numberwithin{theorem}{section}
\numberwithin{prop}{section}
\numberwithin{lemma}{section}
\numberwithin{defin}{section}
\numberwithin{coro}{section}

 \newcommand{\be}{\begin{equation}}
 \newcommand{\ee}{\end{equation}}

\begin{document}

\title[The Coincidence Problem for $\lowercase{d}\leq4$]
{SOLUTION OF THE COINCIDENCE PROBLEM    
              IN DIMENSIONS $\bs{d\leq4}$
\vspace*{1.2cm}}

\author{Michael Baake}
     
\address{Fakult\"at f\"ur Mathematik, Univ.\ Bielefeld, 
                  Pf.\ 100131, 33501 Bielefeld, Germany}
\email{mbaake@math.uni-bielefeld.de}
\urladdr{http://www.math.uni-bielefeld.de/baake/}

\begin{abstract}
  Discrete point sets $\mathcal{S}$ such as lattices or quasiperiodic
  Delone sets may permit, beyond their symmetries, certain isometries
  $R$ such that $\mathcal{S}\cap R\mathcal{S}$ is a subset of
  $\mathcal{S}$ of finite density. These are the so-called {\em
    coincidence isometries}. \index{coincidence isometries} They are
  important in understanding and classifying grain boundaries and
  twins in crystals and quasicrystals.  It is the purpose of this
  contribution to introduce the corresponding coincidence problem in a
  mathematical setting and to demonstrate how it can be solved
  algebraically in dimensions 2, 3 and 4.  Various examples both from
  crystals and quasicrystals are treated explicitly, in particular
  (hyper-)cubic lattices and quasicrystals with non-crystallographic
  point groups of type $H_2$, $H_3$ and $H_4$.  We derive
  parametrizations of all linear coincidence isometries, determine the
  corresponding coincidence index (the reciprocal of the density of
  coinciding points, also called $\gS$-factor), and finally
  encapsulate their statistics in suitable Dirichlet series generating
  functions. 
\end{abstract}

\maketitle

\vspace*{2mm}           
\section{Introduction}

The concept of a {\em \underline{c}oincidence \underline{s}ite
  \underline{l}attice\/} (CSL) arises in the crystallography of grain
and twin boundaries \cite{Boll70,Ranga66,Sass}.  Different domains of
a crystal across a boundary are related by having a sublattice (of
full rank) in common. This is the CSL.  It can be viewed as the
intersection of a lattice with a rotated copy of itself, where the
points in common form a sublattice of {\em finite} index (we shall not
discuss any situation other than that).  Until recently, CSLs have
been investigated only for true lattices or for crystallographic
packings, for example cubic or hexagonal crystals
\cite{Ranga66,Grimmer74,GBW}.  Although the subject itself is quite
old, no systematic investigation in more mathematical terms has been
carried out so far.

With the advent of quasicrystals, many new cases arose, since
quasicrystals also have grain boundaries and one would like to know
the coincidence site quasilattices
\cite{Ranga90,Ranga93,Warrington,WaLu} and, more specifically, which
of them can form twins (or multiple twins, where the angle between the
grains is a rational multiple of $\pi$).  Added impetus is given by
the experimental progress made in recent years \cite{Urban}, in
particular on the study of twins in icosahedral structures.  So, an
extension of the CSL analysis to {\em all}\/ discrete structures is
desirable.

It is the aim of this article to provide the mathematical basis for
it, and to display various relevant examples in detail. Though we
shall include many proofs, it is not possible to present a complete
account here, and we have to refer to original sources several times.
Additional material can be found in \cite{BP94,Pleasants,B96}. Let us
briefly outline how this article is organized.

In Section 2, we set the scene for the periodic case and derive
various results on the coincidence structure of lattices. As we shall
need a generalization to modules later on, most results are formulated
to match that.

Section 3 deals with cubic lattices, where the cases of dimensions 2, 3
and 4 are treated explicitly.  The groups of coincidence isometries
are derived together with the corresponding index formul{\ae}, and the
statistics of the CSLs is encapsulated in Dirichlet series generating
functions which are related to various Dedekind zeta functions.

Sections 4 and 5 consider quasicrystals. Here, the problem must be
split into two parts, one being the coincidence problem for the
underlying limit translation module (which is universal and discussed
in detail) and the other being a problem specific correction in the
projection formalism (which is only briefly outlined).  This is
followed by some concluding remarks.

Finally, in the Appendix, a closely related problem is presented.
While CSLs are special sublattices which depend on metric properties,
the number of {\em all}\/ sublattices of a given index is an affine
property and depends only on the rank. The solution is derived
explicitly for arbitrary rank.

\section{Preliminaries and some general results} \label{sec2}

This paragraph focuses on lattices, although many properties could
directly be formulated for modules.  We prefer this approach, as not
every reader might directly want to go beyond the lattice situation.
We will then generalize the concepts briefly when we pass on to
quasicrystals.

The first concept we need is that of a {\em lattice\/} in
$d$-dimensional Euclidean space (or $d$-space for short), where we follow
the standard approach, see \cite{Cassels} or \cite{Lekker} for details.
\begin{defin}
   A discrete subset $\gG$ of $\EE{}^d$ is called a
   {\em lattice}\/ (of rank and dimension $d$) if it is the
   $\ZZ$-span of $d$ vectors $\bs{a}^{}_1,\ldots,\bs{a}^{}_d$
   that are linearly independent over $\RR$.
   These vectors form a {\em basis}\/ of the lattice.
\end{defin} \index{lattice}
In particular, we can write 
$\gG = \ZZ\bs{a}^{}_1 \oplus \cdots \oplus \ZZ\bs{a}^{}_d$,
and $\gG$ is isomorphic to the free Abelian group of rank $d$.
Another rather common (and equivalent) characterization is to
say that $\gG$ is a co-compact discrete subgroup of $\RR^d$.

Beyond a lattice $\gG$, we shall also need its {\em dual}, 
$\gG^*$, which is given by
\begin{defin}  \label{dual}
   $\;\; \gG^* \, := \, \{ \bs{x} 
          \mid \bs{x} \!\cdot\! \bs{y} \in \ZZ \mbox{ for all } 
          \bs{y} \in \gG \}$.
\end{defin} \index{dual lattice}
Here, $\bs{x} \!\cdot\! \bs{y}$ denotes the standard Euclidean scalar 
product. The dual is of course also a lattice, since 
$\bs{a}_1^*, \ldots , \bs{a}_d^*$
(defined through $\bs{a}_k^* \cdot \bs{a}_{\ell}^{} = \delta_{k\ell}$)
is a basis of $\gG^*$, called the {\em dual basis}.
It is convenient to attach a basis matrix $B$ to a lattice $\gG$,
where the $k$-th column of $B$ consist of the coordinates
of $\bs{a}_k$ in the standard Euclidean basis
$\bs{e}^{}_1, \ldots ,\bs{e}^{}_d$.
It is clear from Definition~\ref{dual} that, if $B$ is a basis
matrix for $\gG$, then $B^*:=(B^{-1})^t$ is one for $\gG^*$.

To proceed, we need the concept of a {\em sublattice}.
Here, in view of generalizations to come in later chapters, 
we employ the group structure.
\begin{defin}
   Let $\gG$ be a lattice in $\EE{}^d$. A subset
   $\gG' \subset \gG$ is called a {\em sublattice}\/ of 
   $\gG$ when it is a subgroup of finite index. The latter
   is the number of residue classes $($or cosets\/$)$ of 
   $\gG'$ in $\gG$, denoted by\/ $[\gG : \gG']$.
\end{defin} 
By definition, a sublattice has finite index in its (super-)lattice,
and thus shares rank and dimension with it.
In particular, $\ZZ^2$ is not a sublattice of $\ZZ^3$ in this notation.
If such a situation occurs, we shall speak of {\em net planes}\/ or
lattice planes (as $\ZZ^2 = \ZZ^3 \cap \{ z=0 \}$).

{}For lattices, the index also has a direct geometric meaning:\ it is
nothing but the (inverse) quotient of the volumes of the fundamental 
domains of the two lattices, as can easily be seen with explicit bases 
for them. A well known result \cite{Cassels} in this context is the
following.
\begin{lemma} \label{sub1}
    Let $\gG$ be a lattice in\/ $\EE{}^d$ with basis matrix $B$.
    Then, $\gG'$ is a sublattice of $\gG$ if and only if
    there exists a non-singular integral matrix $Z$ such that
    $B'=BZ$ is a basis matrix for $\gG'$. 
    The corresponding index is\/ $[\gG : \gG'] = \lvert\det(Z)\rvert$. 
    \qed
\end{lemma}
This description of sublattices is indeed helpful as one can now show
\begin{lemma} \label{subsuper}
  Let $\gG'$ be a sublattice of $\gG\subset\EE^d$ of index $m$.
  Then, $m\gG$ is a sublattice of $\gG'$ of index\/ $m^{d-1}$.
\end{lemma}

\noindent {\sc Proof}: Note that $\gG'\subset\gG$ means $B'=BZ$ with
$Z$ integral and $\lvert\det(Z)\rvert=m$, by Lemma~\ref{sub1}. But
$mB=m B' Z^{-1}=B' Z'$ is a basis matrix for $m\gG$, and $Z' = m
Z^{-1}$ is integral (by the standard formula for the inverse of a
matrix). This means $m\gG\subset\gG'$ by Lemma~\ref{sub1}.  Finally,
$\det(mB)=m^d\cdot \det(B)$ gives the statement about the indices.
\qed

\smallskip Note that, for $d>1$, $m\gG$ need not be the maximal
sublattice of $\gG'$ which is a homothetic copy of $\gG$, as can be
seen from the example $\gG=\ZZ^2$, $\gG'=2\ZZ^2$.  Concerning the dual
lattice, one can show another useful result.
\begin{lemma} \label{sub3}
    If $\gG_2$ is a sublattice of $\gG_1$, one has
    $[ \gG_2^* \, : \, \gG_1^* ] \, = \,
      [  \gG_1 \, : \, \gG_2 ]$, and the
    corresponding factor groups are isomorphic:
    $ \gG_2^*/\gG_1^* \simeq \gG_1^{}/\gG_2^{}$.
\end{lemma}

\noindent {\sc Proof}: 
By assumption, $B_2 = B_1 Z$ with $Z$ integral and non-singular.
Consequently, one has $B_1^* =  B_2^* Z^t$ and the first statement
follows from $\det(Z^t)=\det(Z)$.  The two Abelian factor groups are
thus of equal order.  Isomorphism follows from the observation that
the subgroups of $\gG_1$ which contain $\gG_2$ are, by duality and
part one of the Lemma, in one-to-one relation with the subgroups of
$\gG_2^*$ which contain $\gG_1^*$, with preservation of the
corresponding indices.  \qed

\smallskip
Having prepared the ground, we can now introduce various concepts 
and results which will be helpful for the coincidence problem.
\begin{defin}
   Two lattices $\gG_1,\gG_2$ are called {\em commensurate},
   denoted by $\gG_1 \sim \gG_2$, when $\gG_1\cap\gG_2$
   is a sublattice $($of finite index\/$)$ of both $\gG_1$ and $\gG_2$.
\end{defin}
It is clear that only lattices of the same rank can be commensurate
to one another, and one also has
\begin{prop}
   Commensurateness of lattices is an equivalence relation.
\end{prop}

\noindent {\sc Proof}: 
Reflexivity and symmetry are clear by definition. 
Finally, $\gG_1 \sim \gG_2$ and $\gG_2 \sim \gG_3$  together 
imply that $\gG_1\cap\gG_2\cap\gG_3$ is of finite index in both
$\gG_1$ and $\gG_3$ (one way to see this is via suitable use of
Lemma~\ref{subsuper}), which implies $\gG_1 \sim \gG_3$ 
and hence transitivity.     \qed

\begin{prop} \label{comm2}
    Let $\gG_1$ and $\gG_2$ be commensurate lattices.
    Then, $\gG_1\cap\gG_2$ and $\gG_1 + \gG_2$
    are lattices, and one has the following diagram,
  \begin{displaymath}
    \begin{array}{ccc}
          & \gG_1 + \gG_2 &  \\
          m_1^{}\nearrow & & \nwarrow m_2^{}        \\
          & & \\
          \gG_1  & & \gG_2   \\ 
          & & \\
          n_1^{}\nwarrow & & \nearrow n_2^{}        \\
          & \gG_1 \cap \gG_2 & 
     \end{array}
  \end{displaymath}
 where $A\subset B$ is written as $A\to B$ and the 
 indices satisfy $m_1^{} = n_2^{}$ and $m_2^{} = n_1^{}$.

{}Furthermore, the following equations hold:
  \begin{displaymath}
     (\gG^{}_1\cap\gG^{}_2)^* = \gG_1^* + \gG_2^* 
     \quad \mbox{and} \quad
     (\gG^{}_1 + \gG^{}_2)^* = \gG_1^* \cap \gG_2^* \,.
  \end{displaymath}
\end{prop}

\noindent {\sc Proof}: 
Recall that $\gG_1 + \gG_2$ is defined as
\[
    \gG_1+\gG_2 \; := \; \{ \bs{x}^{}_1 + \bs{x}^{}_2  
    \mid \bs{x}^{}_1\in\gG_1 \, , \; \bs{x}^{}_2\in\gG_2 \} \, .
\] 
It is the smallest group which contains the lattices $\gG_1$ and 
$\gG_2$. As $\gG_1\sim\gG_2$, Lemma~\ref{subsuper} 
guarantees that $\gG_1 + \gG_2$ is again a lattice because
$k (\gG_1+\gG_2)\subset\gG_1\cap\gG_2$ for some $k\in\NN$.

The claim about the indices is a direct consequence of the second 
isomorphism theorem for groups \cite{Zass} giving 
\[
    \gG_2/ (\gG_1\cap\gG_2) \; \simeq \; 
     (\gG_1 + \gG_2)/\gG_1 \quad \mbox{and} \quad
    \gG_1/ (\gG_1\cap\gG_2) \; \simeq \;
     (\gG_1 + \gG_2)/\gG_2 \, .
\] 
As $\gG_1\cap\gG_2$ is of finite index both in $\gG_1$ and $\gG_2$,
all these factor groups are finite.

The second claim follows directly from the definition of the dual 
lattice: $\bs{x}\in(\gG_1+\gG_2)^*$ means 
$\bs{x}\cdot\bs{y}\in\ZZ$ for all $\bs{y}\in(\gG_1+\gG_2)$,
hence in particular for all $\bs{y}\in\gG_k$, $k=1,2$, and
therefore $\bs{x}\in(\gG^*_1\cap\gG^*_2)$. Conversely, if 
$\bs{x}\cdot\bs{y}^{}_k\in\ZZ$
for $\bs{y}^{}_k\in\gG^{}_k$, $k\in\{1,2\}$, then 
$\bs{x}\cdot(\bs{y}^{}_1+\bs{y}^{}_2)\in\ZZ$ and
$(\gG^*_1\cap\gG^*_2)\subset(\gG_1+\gG_2)^*$. 
The other identity follows by duality.     \qed

\smallskip
Now we come to the central concept of this article.
\begin{defin} \label{centraldef}
   Let $\gG$ be a lattice in $\EE{}^d$. An orthogonal
   transformation $R\in \mathrm{O}(d):=\mathrm{O} (d,\RR)$ 
   is called a {\em coincidence isometry}\/
   of $\gG$ when $R\gG \sim \gG$. The integer
   $\gS(R) := [\gG : (\gG\cap R\gG)]$
   is called the {\em coincidence index}\/ of $R$. If $R$ is not
   a coincidence isometry, $\gS(R):=\infty$. We set
\begin{displaymath}
   \begin{array}{ccc}
     \mathrm{OC}(\gG)  &  :=  & \{ R \in \mathrm{O}(d) 
     \mid \gS(R) < \infty \} \\[2mm]
     \mathrm{SOC}(\gG)  &  :=  & 
     \{ R \in \mathrm{OC}(\gG) \mid \det(R)=1 \}
   \end{array}
\end{displaymath}
\end{defin}
If necessary, we shall use the lattice as a subscript for the
coincidence index, $\gS=\gS_{\gG}$, but usually we can safely manage
without.  In this article, we consider only {\em linear\/} isometries,
called isometries for simplicity. In more general situations,
extensions to affine transformations are necessary, compare
\cite[App.~A]{Pleasants}, but we shall not use them here.  One
immediate consequence of Definition~\ref{centraldef} is
\begin{theorem}
   If $\gG$ is a lattice in\/ $\EE^d$,
   $\mathrm{OC}(\gG)$ and $\mathrm{SOC}(\gG)$ are subgroups 
   of\/ $\mathrm{O}(d)$.
\end{theorem}

\noindent {\sc Proof}: 
Let $R_1$ and $R_2$ be in $\mathrm{OC}(\gG)$, i.e., $\gG\sim R_1\gG$ and
$\gG\sim R_2\gG$. Clearly, this also implies $\gG\sim R^{-1}_{2}\gG$ and
$R^{}_{1}\gG\sim R^{}_{1} R^{-1}_{2}\gG$.  From transitivity, we may 
conclude that $\gG\sim R^{}_{1} R^{-1}_{2}\gG$ and hence $R^{}_{1}
R^{-1}_{2}\in\mathrm{OC} (\gG)$.  So, $\mathrm{OC}(\gG)$ is a subgroup of
$\mathrm{O} (d)$; the corresponding statement for $\mathrm{SOC}(\gG)$ is
obvious.  \qed

\smallskip
Without further information on $\gG$, one cannot determine
the corresponding index $\gS=\gS_{\gG}$, except that we know:
\begin{lemma} \label{invindex}
   Let $R\in\mathrm{OC} (\gG)$. Then, $R^{-1}\in\mathrm{OC} (\gG)$
   and $\gS (R) = \gS (R^{-1})$.
\end{lemma}
\noindent {\sc Proof}:
The first claim is clear from the group property of
$\mathrm{OC} (\gG)$. Since $R$ is an isometry, its action
does not change the volume of fundamental domains, and we
obtain the equation
\[
    \gS(R) \; = \;  [\gG:(\gG\cap R\gG)] \; = \; 
    [R\gG:(\gG\cap R\gG)]  \; = \;
    [\gG:(R^{-1}\gG\cap\gG)] \; = \; \gS(R^{-1}).
\]
This establishes the claim. \qed

\smallskip
More can be said about the $\mathrm{OC}$-groups of
related lattices. The following is immediate:
\begin{lemma} 
     If $\gG$ is a lattice, $R$ an orthogonal
     transformation and $\lambda\in\RR\backslash\{0\}$, one has
     the relations\/ $\mathrm{OC}(\lambda\gG) = \mathrm{OC}(\gG)$
     and\/ $\mathrm{OC}(R\gG) = R\,  \mathrm{OC}(\gG)R^{-1}
     \simeq \mathrm{OC}(\gG)$.   \qed
\end{lemma}
Since orthogonal transformations and homotheties commute, the
$\mathrm{OC}$-group does not change (up to conjugation) if one applies
linear similarity transformations.  This is one step to show that the
$\mathrm{OC}$-group essentially is an invariant of so-called Bravais
classes (cf.\ \cite{Schwarz} for details on this concept) -- if one
disregards non-generic solutions for special representatives. The
latter problem does not occur if one deals with irreducible
symmetries, such as fourfold symmetry in the plane or full cubic
symmetry in $\EE^3$. In this article, no other situation shall be
discussed.

\begin{lemma}
   Let $\gG_2$ be a sublattice of $\gG_1$ of 
   $($finite\/$)$ index\/ $m = [\gG_1 : \gG_2]$.
   Then, $\mathrm{OC}(\gG_2) = \mathrm{OC}(\gG_1)$ and,
   for any isometry $R$ out of this group, one has\/ 
   $\gS_1 (R) \, | \; m  \gS_2 (R)$.
\end{lemma}

\noindent {\sc Proof}: 
If $R$ is a coincidence isometry of $\gG_2$, we know that
$\gG_2\cap R\gG_2$ has finite index in $\gG_2$ 
and hence also in $\gG_1$. As
$\gG_2\cap R\gG_2 \subset \gG_1\cap R\gG_1 
\subset \gG_1$, 
we see $\mathrm{OC}(\gG_2) \subset \mathrm{OC}(\gG_1)$.
On the other hand, $[\gG_1:\gG_2]=m$ implies 
$m\gG_1\subset\gG_2$ by Lemma~\ref{subsuper}, so also 
$\mathrm{OC}(\gG_1) \subset \mathrm{OC}(\gG_2)$
and the two groups must be equal.
The second statement follows from 
$\gG_2\cap R\gG_2 \subset 
\gG_1\cap R\gG_1 \subset \gG_1$ and
$\gG_2\cap R\gG_2 \subset \gG_2
\subset \gG_1$ by a comparison of indices.   \qed

\smallskip
If we write $\gG_1 = \dot{\bigcup}_{\ell=0}^{m-1} (t_{\ell}
+\gG_2)$ with $t_0=0$, one can say more about the relation
between $\gS_1 (R)$ and $\gS_2 (R)$ if $R$ respects the cosets
in suitable ways. In particular, if $R\gG_2$ is disjoint from
the cosets $t_{\ell}+\gG_2$ for all $\ell>0$, one can also
derive that $\gS_2 (R) | \gS_1 (R)$.

The last results on the coincidence groups can obviously be extended to
\begin{coro}
   Commensurate lattices possess the same $\mathrm{OC}$-group.  \qed
\end{coro}

With the above results, one can now also relate a lattice with its dual.
\begin{theorem} \label{sigdual}
   Let $\gG$ be a lattice in\/ $\EE^d$ and $\gG^*$ its dual. Then,
   $\mathrm{OC}(\gG^*) = \mathrm{OC}(\gG)$ and the coincidence index
   of any orthogonal matrix is the same for both lattices.
\end{theorem}

\noindent {\sc Proof}: 
Assume $\gG \sim R\gG$. Then, since
$\gG^{**} = \gG$ and $(R\gG)^* = R \gG^*$,
we have $\infty > m = \gS(R)
= [\gG : (\gG\cap R\gG)] 
= [(\gG + R\gG) : \gG]
= [\gG^* : (\gG + R\gG)^*] 
= [\gG^* : (\gG^* \cap R\gG^*)] 
= \gS^*(R)$, by applying Proposition 
\ref{comm2} (specialized to the case $\gG_2=R\gG_1$,
where $m^{}_1=m^{}_2=n^{}_1=n^{}_2$) and Lemma~\ref{sub3}.
It also follows that $\gS(R)=\infty$ if and
only if $\gS^*(R)=\infty$.  \qed

\smallskip We have now all prerequisites to tackle the coincidence
problem for lattices and crystals.  While we proceed, we shall need,
in particular in the part on quasicrystals, various results from
algebraic number theory. It is not possible to present a
self-contained description here, but references to the relevant
literature shall be given. For general background, we refer to
\cite[Ch.~1.4]{Lekker}, to \cite[Chs.\ XIV--XVII and Ch.\ XX]{Hardy},
and, for a number theoretic approach to quasicrystals in general, to
\cite{PABP}.

\section{Lattices and crystals: the cubic case}

Although there are several cubic lattices, let us first consider
the primitive (hyper-)cubic lattice, $\ZZ^d$. With the standard
Euclidean basis $\bs{e}^{}_{1},\ldots,\bs{e}^{}_{d}$, it can be
written as 
\[
    \ZZ^d  \; = \; 
   \ZZ\bs{e}^{}_{1}\oplus\dots\oplus\ZZ\bs{e}^{}_{d}.
\]
As we shall need the term `primitive' later on in a different meaning,
we replace it here by the term $P$-type from now on, as compared
to $F$-type (for face centred cubic, {\it fcc}, and its generalizations)
and to $B$-type (for body centred, {\it bcc}, and its generalizations).

It is immediately clear that an orthogonal matrix $R$ with rational entries
only is a coincidence isometry of $\ZZ^d$: one can directly give a lattice
which is a common sublattice to both $\ZZ^d$ and $R\ZZ^d$, namely $m \ZZ^d$,
where $m$ is the {\em denominator}\/ of $R$ defined through 
\be \label{den}
   \den(R) \; := \; \gcd\{k \mid k \cdot R 
     \mbox{ integral} \} \, ,  
\ee 
where $\gcd$ is the greatest common divisor. On the
other hand, as soon as one entry of $R$ is {\em irrational}, $R_{ij}$ say, no
point in the direction of $\bs{e}_j$ coincides with a lattice point after
rotation. This means that we have one lattice direction without any
coincidence, hence infinitely many residue classes of the set of coinciding
points which therefore is not of finite density, and $R$ cannot be a
coincidence isometry. So we have 
\begin{theorem} \label{cubic1} $\;\;$
      $ \mathrm{OC}(\ZZ^d) \, = \,  \mathrm{O}(d,\QQ)$ and\/ 
      $\mathrm{SOC}(\ZZ^d) \, = \, \mathrm{SO}(d,\QQ)$.   \qed
\end{theorem}

Having this result for the $P$-type (hyper-)cubic lattices, the 
obvious next question is what happens for the other lattices with
(hyper-)cubic symmetry. The types of such lattices, up to similarity
transformations, are as follows. 
\begin{theorem}
     In dimensions $d=3$ and $d\geq5$, there are precisely
     three $($hyper-$)$cubic lattices, namely $F$-type, $P$-type 
     and $B$-type. For $d=1$ and $d=2$, there is
     only one such lattice $($represented by the integers resp.\ by
     the square lattice\/$)$, while in $d=4$, 
     there are two such lattices $(P$-type hypercubic and centred\/$)$.  \qed
\end{theorem}

A proof can be found in \cite{Schwarz}.  It is possible, in any dimension, to
realize the (hyper-) cubic lattices in such a way that they are commensurate
-- usually, one works with $\ZZ^d$ as representative of $P$-type, with the
root (weight) lattice $D^{}_d$ ($D^{*}_{d}$) as representatives of $F$-type
($B$-type) in dimensions $d\ge 5$, with $D^{}_4$ as centred lattice in
$4$-space, and with the usual {\it fcc}\/ and {\it bcc}\/ lattices in
$3$-space (or $A^{}_3$ resp.\ $A^{*}_{3}$), compare \cite{Conway}.
With this convention, the $\mathrm{OC}$-groups of the (hyper-)cubic 
lattices in $\EE^d$ (with $d$ fixed) are
identical, though the corresponding indices might differ.  At present, we do
not know the general answer in higher dimensions, but more can be said about
the cubic lattices in dimensions $2$, $3$ and $4$. So, let us summarize some
of those results.

\subsection{$\bs{d=2}$: the square lattice}

Let us first describe the case of the square lattice in more detail,
as this is the simplest non-trivial example. Here, we shall be less
formal and refer for all the details and proofs to 
\cite{Pleasants} without any further mentioning. First, we restrict
the description to rotations, and come to reflections at the
end of the section.

The square lattice $\ZZ^2$, embedded in $\EE{}^2$ resp.\ $\RR^2$,
consists of all integer linear combinations of the two
vectors $\bs{e}_1$ and $\bs{e}_2$, i.e., $\ZZ^2 =
\ZZ\bs{e}^{}_{1} \oplus \ZZ\bs{e}^{}_{2}$.
A rotated copy $R \, \ZZ^2$ with $R  \in  \mathrm{SO}(2,\RR)$ 
results in a CSL of finite index, according to Theorem~\ref{cubic1}, 
if and only if both $\cos(\varphi)$ and
$\sin(\varphi)$ are {\em rational}, where $\varphi$ is the rotation angle.
This gives the well-known relation between
coincidence rotations and primitive Pythagorean triples \cite{Lueck}. 
The group of coincidence rotations is thus explicitly seen to be
$    \mathrm{SOC}(\ZZ^2)  =  \mathrm{SO}(2,\QQ)  $.
To investigate this group, we employ some elementary results
from the algebraic theory of quadratic fields
\cite[Chs.\ XIV and XV]{Hardy}.
In particular, we notice that, with $i = \sqrt{-1}$, 
we can identify $\ZZ^2$ with the ring of Gaussian integers 
(the algebraic integers of the quadratic field $\QQ(i)$, an
extension of $\QQ$ of degree $2$):
\begin{equation}
   \ZZ^2 \; = \; \ZZ[i] \; = \; \{ m+n i \, | \, m,n \in \ZZ \} \,.
\end{equation}
The ring $\ZZ[i]$ is a Euclidean domain and thus has unique
factorization up to units, see \cite{Hardy}.
The units ($i$ and its powers) form a group isomorphic to $C_4$,
the rotation part of Aut($\ZZ^2$).

In this setting, a rotation $R(\varphi) \in \mathrm{SOC}(\ZZ^2)$
corresponds to multipli\-cation by a complex number $e^{i\varphi} \in
\QQ(i)$.  This number can be written as $e^{i\varphi} = \alpha/\beta$
with $\alpha,\beta \in \ZZ[i]$ coprime and of equal norm, i.e.,
$\lvert \alpha\rvert = \lvert\beta\rvert$ (so, $e^{i\varphi}$ rotates
the lattice point $\beta$ into another lattice point, $\alpha$, on the
same circle around the origin). We now factorize numerator and
denominator into Gaussian primes. They are the rational (or ordinary)
primes $p\equiv 3$ mod $4$, the factor $1+i$ of $2$ (a so-called
ramified prime), and the pairs of complex conjugate factors of
rational primes $p\equiv 1$ mod $4$ (where $p=\omega_p\,
\overline{\omega}_p$, and $\omega_p/\overline{\omega}_p$ is not a unit
in $\ZZ[i]$). Clearly, as $\alpha$ and $\beta$ are coprime and both
divide the same rational integer $\ell=\lvert\alpha\rvert^2=
\lvert\beta\rvert^2$, only the last type of primes can occur in the
factorization, always one (or a power of it) in the numerator and its
complex conjugate in the denominator.

As a consequence of the unique factorization property up to units,
every coincidence rotation can then be factorized as
\begin{equation} \label{product1}
   e^{i\varphi} \, = \, \varepsilon \cdot \!\!\!
                \prod_{p \equiv 1 \; (4)}
         \left( \frac{\omega_p}{\overline{\omega}_p} \right)_{ }^{n_p}
\end{equation}
where $n_p\in\ZZ$ (only finitely many of them $\neq 0$), $\varepsilon$
is a unit in $\ZZ[i]$ (a power of $i$), $p\/$ runs through the
rational (or ordinary) primes congruent to 1 (mod 4), and the
$\omega_p$, $\overline{\omega}_p$ are the (complex conjugate) Gaussian
prime factors of $p$. We thus have the (non-trivial!) result that
$\mathrm{SOC} (\ZZ^2)$ is an infinitely generated Abelian group that
nevertheless permits the factorization into a torsion group and a free
Abelian group, namely
\begin{equation}  \label{dirprod1}
      \mathrm{SOC}(\ZZ^2) \, = \, \mathrm{SO}(2,\QQ)
     \, \simeq \,  C_4 \times \ZZ^{(\aleph_0)}
\end{equation}
with generators $i$ for $C_4$ and $\omega_p / \overline{\omega}_p$
with $p \equiv 1 \; (4)$ for the infinite cyclic groups.
By $\ZZ^{(\aleph_0)}$ we mean, as usual, the infinite Abelian group
that consists of all {\em finite} integer linear combinations in the
(countably many) generators. We prefer a multiplicative rather
than additive notation here as the coincidence groups will not
be Abelian in later examples.

To determine the coincidence index $\gS(R)$, we observe that it is $1$
for the units (i.e., true symmetry rotations) and $p$ for the
generator $\omega_p / \overline{\omega}_p $ (since $p=
\omega_p\,\overline{\omega}_p = \mathrm{norm} (\omega_p)$ counts
the number of residue classes of the corresponding CSL in $\ZZ^2$).
More generally, due to multiplicativity, $\gS (R)$ is the
(number theoretic) norm of the numerator of \eqref{product1}, i.e.,
\begin{equation} \label{sigform1}
    \gS(R) \, = \, \prod_{p \equiv 1 \; (4)}
                             p_{ }^{|n_p|} \; .
\end{equation}
This shows the power of the generator approach in this (Abelian)
situation. The first few generators with $\gS > 1$ are
\begin{displaymath} 
 \frac{4+3i}{5} \; , \; \frac{12+5i}{13} \; , \; 
 \frac{15+8i}{17} \; , \;
 \frac{21+20i}{29} \; , \; \frac{35+12i}{37} \; , \;
 \frac{40+9i}{41} \; , \; \mbox{etc.}
\end{displaymath}
These are normalized (by multiplication with a suitable unit) to
have argument in $(0,\pi/4)$, and are shown with
denominator $\gS$ (a prime $\equiv$ 1 (mod 4)).
All other coincidence rotations are obtained by
combinations, as indicated in Eq.~(\ref{product1}).

It is convenient to summarize the possible coincidence indices and the 
number of CSLs with a given index by means of a generating 
function. To do so,
let $4 f(m)$ denote the number of coincidence rotations of index $m$,
which means that $f(m)$ counts the different CSLs of index $m$.
As a consequence of unique factorization in $\ZZ[i]$,
$f(m)$ is a {\em multiplicative}\/ arithmetic function (i.e., $f(1)=1$ and
$f(m_1m_2)=f(m_1)f(m_2)$ for coprime $m_1,m_2$), and we can calculate
the numbers $f(m)$ if we know them for $m=p^r$ with $p$ prime and
$r>1$. This simplification is a nice algebraic result that need
no longer hold for the analogous problem applied to planar
modules with $N$-fold symmetry when $N\ge 46$, see \cite{Pleasants}.

In our present case, due to multiplicativity, a
Dirichlet series $\varPhi(s)$ is an appropriate generating 
function \cite{Apostol,Wilf}. To calculate $f(m)$ explicitly,
we observe that, if $m=p^r$ is a prime power ($p\equiv 1$ (mod $4$),
$r\ge 1$), only the two rotations
\[
    \left( \frac{\omega_p}{\overline{\omega}_p} \right)^r  \; , \quad
    \left( \frac{\overline{\omega}_p}{\omega_p} \right)^r
\]
lead to CSLs of index $m$, wherefore we have $f(p^r)=2$ in
this case. Now, we can directly calculate the entire generating
function through its Euler product representation and obtain:

\begin{prop} \label{index1}
   The Dirichlet series generating function of the number $f(m)$
   of CSLs of\/ $\ZZ^2$ of index $m$ is
\[
\begin{split}
\varPhi(s) & \; = \; \sum_{m=1}^{\infty}\frac{f(m)}{m^s}
  \; = \; \prod_{p \equiv 1 \, (4)} \Big(1+\frac{2}{p^s}+
                \frac{2}{p^{2s}} +\cdots \Big) 
  \; = \; \prod_{p \equiv 1 \, (4)}\frac{1+p^{-s}}{1-p^{-s}} \\
& \; = \; 1+\tfrac{2}{5^s}+\tfrac{2}{13^s}+\tfrac{2}{17^s}+
          \tfrac{2}{25^s}+\tfrac{2}{29^s}+\tfrac{2}{37^s}+
         \tfrac{2}{41^s}+\tfrac{2}{53^s} + \tfrac{2}{61^s} + 
    \tfrac{4}{65^s} + \tfrac{2}{73^s} + \cdots 
 \end{split}
\]
\end{prop}

A little later, we shall express $\varPhi(s)$ in terms of
$\zeta$-functions.  This generating function is not only a succinct
way of representing the statistics of CSLs and coincidence indices, 
it is also
a powerful tool for determining their asymptotic properties.  For
example, it can be used to show (through the determination of the
right-most pole of $\varPhi (s)$ in the complex $s$-plane, which is at
$s=1$, and its residue) that the number of CSLs of $\ZZ^2$ with index
$\le N\/$ is asymptotically $N/\pi$ (and the corresponding number of
coincidence rotations is asymptotically $4N/\pi$).  The possible
coincidence indices are precisely the numbers $m$ with all prime
factors $\equiv 1 \; (\mbox{mod}\, 4)$, and we then have \be f(m) \; =
\; 2^a \, , \ee where $a$ is the number of distinct prime divisors of
$m$.  Each CSL is itself a square lattice \cite{Pleasants,BG}, with 
the index as the area of its fundamental domain. We shall come back 
to this in a more general context.

{}Finally, the full group of coincidence isometries,
$\mathrm{OC}(\ZZ_{ }^2)$, is the semidirect
product of the rotation part 
$\mathrm{SOC}(\ZZ_{ }^2)$ (normal subgroup) with the cyclic
group $C_2$ generated by complex conjugation ($=$ reflection in the
$x$-axis):
\be \label{isom1}
     \mathrm{OC}(\ZZ_{ }^2) \; = \;  \mathrm{SOC}(\ZZ_{ }^2) \rtimes C_2 \; .
\ee
Here, conjugation of a rotation through an angle $\varphi$ by complex
conjugation results in the inverse rotation (through $-\varphi$).
Let us give a justification of Eq.~(\ref{isom1}).
Since $\mathrm{O}(2) = \mathrm{SO}(2) \rtimes C_2$ (semi-direct product)
with the $C_2$ of Eq.~(\ref{isom1}),
any planar isometry $T$ with $\det(T)=-1$ can uniquely be written
as the product
\be \label{factor1}
     T \, = \, R(\varphi) \cdot T_x \,
\ee
of a rotation through $\varphi$ with $T_x$, the reflection in 
the $x$-axis (complex conjugation).
But $T_x$ leaves the entire lattice $\ZZ_{ }^2$ invariant, 
whence $T$ is a coincidence isometry
if and only if $R(\varphi)$ is a coincidence rotation.

The calculation of coincidence indices is also simple in this case.
The coincidence index for the reflection $T_x$ is 1.
{}For an arbitrary element of $\mathrm{OC}(\ZZ_{ }^2)$, we either 
meet a rotation (where we know the result already) or use the 
factorization of Eq.~\eqref{factor1} again.
Then, the coincidence index is identical with that of its rotation part,
so Eq.~(\ref{factor1}) is all that is needed.   
This solves the coincidence problem for the square lattice completely
and we have
\begin{theorem}
   The group of coincidence isometries of the square lattice $\ZZ^2$ is
   \begin{displaymath}
      \mathrm{OC}(\ZZ^2) \; = \; \mathrm{O}(2,\QQ) \; \simeq \; 
          (C_4 \times \ZZ^{(\aleph_0)}) \, \rtimes \, C_2 \; .
   \end{displaymath}  
   This group is fully characterized by Eqs.~$\eqref{product1}$, 
   $\eqref{dirprod1}$ and $\eqref{isom1}$. The coincidence index of a 
   rotation is given by Eq.~$\eqref{sigform1}$, while that of a reflection 
   $\eqref{factor1}$  equals the index of its rotation part.   
   The corresponding Dirichlet series generating function is\/ 
   $\varPhi (s)$ from Proposition~$\ref{index1}$. \qed
\end{theorem}

\subsection{A short digression on a hierarchy of problems}

It is the intention of
this paragraph to shed some more light on the coincidence
problem and how it relates to similar questions. We shall
explain it for the square lattice in an informal manner.

To this end, let us start with the question of how many sublattices
of $\ZZ^2$ have index $m$ -- without any further restriction.
Let us call this number $a_m$. Clearly, $a_1^{}=1$ (only
$\ZZ^2$ itself is sublattice of index 1) and $a_2^{}=3$ (counting
two different rectangular sublattices and one square sublattice).
In general, $a_{mn} = a_m a_n$ when $m,n$ are coprime,
and one can derive, either from the Appendix or from 
\cite[Lemma~2 on p.\ 99]{Serre}, the general result
\[
    a_m \; = \; \sigma_1^{}(m) \; = \; 
    \mbox{$\sum_{d\mid m}$}\, d 
\]
with generating function
\be
\begin{split}
  F(s) & \; = \; \sum_{m=1}^{\infty} \frac{a_m}{m^s}
  \; = \;  \zeta(s) \cdot \zeta(s-1)  \\[2mm]
  & \; = \;
  1+\tfrac{3}{2^s}+\tfrac{4}{3^s}+\tfrac{7}{4^s}+\tfrac{6}{5^s}+
  \tfrac{12}{6^s}+
  \tfrac{8}{7^s}+\tfrac{15}{8^s}+\tfrac{13}{9^s}+\tfrac{18}{10^s}+
  \tfrac{12}{11^s}+\tfrac{28}{12^s} + \tfrac{14}{13^s} + \cdots
\end{split}
\ee 
Here, $\zeta(s) = \sum_{m=1}^{\infty} m^{-s}$ is Riemann's zeta
function, compare \cite{Hardy}.  From this, it can be shown that the
number of sublattices with index $m$ grows, on average, linearly as $m
\pi^2/6$, compare \cite[Thm.~324]{Hardy}.  These results are, of
course, {\em affine} in nature and apply to any lattice of rank 2, and
to any free Abelian group of rank 2 (counting the different subgroups
of index $m$).

Let us look for {\em metric} properties by asking how many
of the sublattices of $\ZZ^2$ of index $m$ are actually {\em square}
lattices (see \cite{BG2,BM1,BM} for various generalizations of this
question). This number can be obtained by counting
the lattice points on circles of radius $m$ (hence
counting solutions of the Diophantine equation $x^2+y^2=m$)
and afterwards dividing by 4 (the order of $C_4$, the rotation
part of the point symmetry group of $\ZZ^2$). The result is
given in Chapters 16.9, 16.10 and 17.9 of \cite{Hardy}
and leads to the Dirichlet series generating function
\be
\begin{split}
  F(s) & \; = \; \zeta^{}_K(s)  \; = \;
            \frac{1}{1-2^{-s}} \cdot
            \prod_{p \equiv 1 \; (4)} 
                   \left(\frac{1}{1-p^{-s}}\right)^2 \cdot
            \prod_{p \equiv 3 \; (4)} \frac{1}{1-p^{-2s}}  \\[2mm]
& \; = \; 1+\tfrac{1}{2^s}+\tfrac{1}{4^s}+\tfrac{2}{5^s}+\tfrac{1}{8^s}+
  \tfrac{1}{9^s}+\tfrac{2}{10^s}+\tfrac{2}{13^s}+\tfrac{1}{16^s}+
  \tfrac{2}{17^s}+\tfrac{1}{18^s}+
  \tfrac{2}{20^s}+\tfrac{3}{25^s} + \cdots 
\end{split}  
\ee 
where here and in what follows $\zeta^{}_K(s)$ is the Dedekind
zeta function of the quadratic field $K=\QQ(i)$, compare \cite[$\S$
63, A.~14 on p.~251]{Scheja}.  The average value of the coefficients
of $F(s)$ is constant, namely $\pi/4$, which follows either from the
asymptotic properties of the generating function near its right-most
pole (at $s=1$) or just from counting one quarter of the lattice
points inside the circle of radius $\sqrt{m}$, which is $m\pi/4$ to
leading order in $m$.

In the last case, some of the square lattices fail to be {\em primitive}
(such as $3 \ZZ^2$ etc.), i.e., whenever the sublattice is an integer
multiple of $\ZZ^2$ or one of its primitve sublattices. If we exclude 
those, primes $p\equiv3 \; (4)$ are impossible as divisors of the index
$m$, and some solutions of $p\equiv1 \; (4)$ also
drop out (whenever the index is divisible by a square). 
Now, the generating function reads
\be
\begin{split}
  F(s) & \; = \; \left( 1 + 2^{-s}\right) \cdot
                 \prod_{p\equiv1 \; (4)}
                 \frac{1+p^{-s}}{1-p^{-s}}
                 \; = \;  \frac{\zeta^{}_K(s)}{\zeta(2s)}  \\[2mm]
& \; = \; 1+\tfrac{1}{2^s}+\tfrac{2}{5^s}+\tfrac{2}{10^s}+
  \tfrac{2}{13^s}+\tfrac{2}{17^s}+
  \tfrac{2}{25^s}+\tfrac{2}{26^s}+\tfrac{2}{29^s}+\tfrac{2}{34^s}+
  \tfrac{2}{37^s}+\tfrac{2}{41^s} + \cdots
\end{split}  
\ee
and the average number of primitive square sublattices
of index $m$ is given by $3/(2\pi)$. This can be determined
by counting one quarter of the visible points \cite{Apostol} in the
circle of radius $\sqrt{m}$, which is $(\pi m \cdot 6/\pi^2)/4$.

Still, not all primitive square sublattices are CSLs -- they only are if the
index is {\em odd}. This finally results in the generating function of the
coincidence problem described above, namely 
\be 
\varPhi(s) \; = \;
\prod_{p\equiv1 \; (4)} \frac{1+p^{-s}}{1-p^{-s}} \; = \; \left( 1 +
  2^{-s}\right)^{-1} \frac{\zeta^{}_K(s)}{\zeta(2s)} 
\ee 
where the average number of CSLs of index $m$ is $1/\pi$.  This equation
expresses the Dirichlet series generating function in terms of zeta functions.
Similar formulas will also appear in later examples.  We hope that this short
digression has put the problem in a broader perspective. Let us now climb up
to higher dimension where the picture changes significantly because
$\mathrm{O}(d)$ is no longer Abelian for $d>2$.

\subsection{$\bs{d=3}$: the three cubic lattices}

In this paragraph, we shall use the notation $\gG_{F,P,B}$ for
the $F$-type ({\it fcc}), $P$-type, and $B$-type ({\it bcc}) cubic 
lattices, respectively. In particular, $\gG_P = \ZZ^3$.
Let us start with this and write the lattice as
$\ZZ^3 = \ZZ \bs{e}_1 \, \oplus \, \ZZ \bs{e}_2 \, 
\oplus \, \ZZ \bs{e}_3$. 
We already know that $\mathrm{OC}(\ZZ^3) = \mathrm{O}(3,\QQ) $
(other isometries might also lead to coincidences, but not
to a CSL of full rank 3).
The subgroup of rotations with index 1 
is the rotation symmetry group of the cube  
of order 24, $\mathrm{Aut}^{+} (\ZZ^3)=
\mathcal{O}=\mathrm{SO}(3,\ZZ)$, cf.\  \cite{Baake84}
for details.

As is well known, $\mathrm{O}(3) = \mathrm{SO}(3) \times C_2$ is a
direct product, where $C_2 = \{\pm\mathbbm{1}_3\}$ is 
the centre of $\mathrm{O}(3)$.
Consequently, we can restrict our attention to pure rotations, as
reflections may be written as a product of a rotation $R$ with
$-\mathbbm{1}^{}_{3}$ (note that $-\ZZ^3 = \ZZ^3$). At this point, we
introduce quaternions $\bs{q}$, see \cite{Hurwitz,duVal,Hardy}, and
Cayley's parametrization \cite{Koecher} with $4$ real numbers
$(\kappa,\lambda,\mu,\nu) = \bs{q} \neq \bs{0}$, 
\be 
   R(\bs{q}) \, = \,
   \frac{1}{|\bs{q}|^2} \mbox{\small $ \left( \begin{array}{ccc} \kappa^2
      + \lambda^2 - \mu^2 - \nu^2 &
      -2\kappa\nu + 2\lambda \mu & 2\kappa\mu + 2\lambda \nu \\
      2\kappa\nu + 2\lambda \mu & \kappa^2 - \lambda^2 + \mu^2 - \nu^2
      &
      -2\kappa\lambda + 2\mu\nu \\
      -2\kappa\mu + 2\lambda \nu & 2\kappa\lambda + 2\mu\nu & \kappa^2
      - \lambda^2 - \mu^2 + \nu^2
    \end{array} \right) $ },
\ee
where $|\bs{q}|^2 = \kappa^2 + \lambda^2 + \mu^2 + \nu^2$ is the
so-called {\em reduced norm}\/ of the quaternion $\bs{q}$. 
In particular, the multiplicative
group of quaternions of norm $1$, which form the unit sphere 
$\mathbb{S}^3$, provides
the usual double cover of the rotation group $\mathrm{SO}(3,\RR)$,
via the group homomorphism $\bs{q}\mapsto R(\bs{q})$
(since $R(\bs{q}) = R(-\bs{q})$). 

As we are interested in rotation matrices with rational entries (i.e., in
the subgroup $\mathrm{SO}(3,\QQ)$), we consider Cayley's parametrization
with $4$ {\em integers}\/ $(\kappa,\lambda,\mu,\nu)=\bs{q}\neq \bs{0}$, and
choose them {\em coprime}, i.e., $\gcd(\kappa,\lambda,\mu,\nu)=1$.
We call such quaternions {\em primitive}. This way, we parametrize the 
entire group $\mathrm{SO}(3,\QQ)$, and obtain each element of it exactly
twice (again because of $R(\bs{q}) = R(-\bs{q})$). With this approach,
one obtains the following result.
 
\begin{prop} \label{cubind}
   Consider $\gG^{}_{P}=\ZZ^3$ and
   let $R=R(\bs{q})\in \mathrm{SO}(3,\QQ)$ be parametrized by a primitive
   quaternion $\bs{q}$. Its coincidence index $\gS(R)$ is the 
   denominator of $R$ as defined in Eq.~$\eqref{den}$. 
   It is the ``odd part" of\/ $|\bs{q}|^2$,
\[
     \gS(R) \, = \, \mathrm{den}(R(\bs{q})) \, = \,
               |\bs{q}|^2/2^{\ell},
\]
where $\ell$ is the largest integer such that\/ $2^{\ell}$ divides\/
$|\bs{q}|^2$.
\end{prop}

Before we prove the statement, let us remark that the necessity for
the division by $2^{\ell}$ in Proposition~\ref{cubind} stems from the
fact that some primitive quaternions (such as $(1,1,0,0)$ or
$(1,1,1,1)$) have norms that are not coprime with the matrix entries
of $R(\bs{q})$. A closer inspection, using results of \cite{Hurwitz},
shows that only the exponents $\ell\in\{0,1,2\}$ are possible.

\smallskip
\noindent {\sc Proof}:
Consider $R(\bs{q})$ with $\bs{0}\neq\bs{q}=(\kappa,\lambda,\mu,\nu)$ 
primitive and $\sigma = \lvert\bs{q}\rvert^2$. Define the $4$ vectors
\[
  \bs{v}^{}_{0} = \begin{pmatrix} \lambda \\ \mu \\ 
  \nu \end{pmatrix}\, , \quad
  \bs{v}^{}_{1} = \begin{pmatrix} \kappa \\ \nu \\ 
  -\mu \end{pmatrix}\, , \quad
  \bs{v}^{}_{2} = \begin{pmatrix} -\nu \\ \kappa \\ 
  \lambda \end{pmatrix}\, ,\quad
  \bs{v}^{}_{3} = \begin{pmatrix} \mu \\ -\lambda \\ \kappa \end{pmatrix}\, ,
\]
which all have integer preimages under $R=R(\bs{q})$, and the $4$ matrices
\[ 
\begin{split}
     B^{}_{0} = (\bs{v}^{}_{1},\bs{v}^{}_{2},\bs{v}^{}_{3}) \, , \quad
     & \det(B^{}_{0}) = \kappa \sigma\, ; \qquad
     B^{}_{1} = (\bs{v}^{}_{0},\bs{v}^{}_{2},\bs{v}^{}_{3}) \, , \quad
     \det(B^{}_{1}) = \lambda \sigma\, ; \\
     B^{}_{2} = (\bs{v}^{}_{1},\bs{v}^{}_{0},\bs{v}^{}_{3}) \, , \quad
     & \det(B^{}_{2}) = \mu \sigma\, ; \qquad
     B^{}_{3} = (\bs{v}^{}_{1},\bs{v}^{}_{2},\bs{v}^{}_{0}) \, , \quad
     \det(B^{}_{3}) = \nu \sigma\, .
\end{split}
\]
Now, each $B_j$ can be read as a basis matrix of a lattice that is a
sublattice of both $\ZZ^3$ and $R \ZZ^3$, hence it is also
a sublattice of the CSL $\ZZ^3\cap R \ZZ^3$.
Consequently, the coincidence index $\gS(R)$ must divide each
of the determinants $\det(B_j)$. As $\bs{q}$ was primitive,
$\gS(R)$ must therefore divide $\sigma$. Since $\gS(R)$ trivially
also divides the third power of $\mathrm{den}(R)=\sigma/2^\ell$
(which is odd), we have $\gS(R) | \mathrm{den}(R)$.


To establish the claim, it is now sufficient to show that also
$\mathrm{den} (R) | \gS (R)$. Observe that
$[\ZZ^3 : \ZZ^3\cap R \ZZ^3] = [ R\ZZ^3 : \ZZ^3\cap R\ZZ^3]$
due to $\det (R)=1$, as in the proof
of Lemma~\ref{invindex}. Since, by definition of the denominator,
\[
   \gcd \{ \mathrm{den} (R) R_{ij}\mid 1\le i,j \le 3\} 
   \; = \; 1 \, ,
\]
the lattice $R\ZZ^3$ contains a vector of the form
$\bs{a}/\mathrm{den} (R)$ with a primitive vector $\bs{a}\in\ZZ^3$,
i.e., $\bs{a}=(a^{}_1, a^{}_2, a^{}_3)^t$ with $a_i\in\ZZ$ and and
$\gcd(a^{}_1, a^{}_2, a^{}_3)=1$. This implies that the number of cosets
of $\ZZ^3\cap R\ZZ^3$ in $R\ZZ^3$ must be a multiple of $\mathrm{den}
(R)$, so that $\mathrm{den} (R) | \gS (R)$.  \qed

\smallskip
In particular, this reproduces the well-known result 
\cite{Grimmer73,Grimmer74} 
that $\gS(R) = \gS_{\ZZ^3} (R)$ runs precisely through all odd integers, i.e.,
$\gS_{\ZZ^3} (\mathrm{SO} (3,\QQ))= 2\NN^{}_{0} + 1 = \{1,3,5,7,\ldots\}$.

Cayley's parametrization has the nice property that 
$\bs{v}^{}_{0} = (\lambda,\mu,\nu)^t$ gives
the (generic) rotation axis of $R(\bs{q}) = R(\kappa,\lambda,\mu,\nu)$,
\begin{equation}
      R(\kappa,\lambda,\mu,\nu) \, \bs{v}^{}_{0}
          \; = \; \bs{v}^{}_{0} \, ,
\end{equation}
while the rotation angle follows from
$\mbox{tr}(R) = 1 + 2\cos(\varphi)$, so
\begin{equation}
  \cos(\varphi) \; = \;
   \frac{\kappa^2 - \lambda^2 - \mu^2 - \nu^2}
        {\kappa^2 + \lambda^2 + \mu^2 + \nu^2} \; .
\end{equation}
One can easily construct all solutions for small indices explicitly,
while the case $\kappa=0$ gives all coincidence rotations through
$\pi$ as described by L\"uck \cite{Lueck}.  If $24f(m)$ is the number
of coincidence rotations of index $m$, the arithmetic function $f(m)$ 
once again counts the
different CSLs of index $m$.  This function is multiplicative as a
consequence of the fact that integer quaternions\footnote{They
  constitute the Hurwitz ring $\mathbb{J}$ of quaternions of the form
  $\frac{1}{2}(a,b,c,d)$ with $a,b,c,d$ all even or all odd.}  have
unique left (and right) factorization up to units\footnote{They are the $24$
  quaternions $\pm (1,0,0,0)$ (and permutations) and $\frac{1}{2}(\pm
  1,\pm 1,\pm 1,\pm 1)$.  Together, they form a group that is the
  double cover of the symmetry group (rotations only) of the regular
  tetrahedron, see \cite{BM}.}, see \cite{Z2,BLP96}
for details.  Its calculation amounts to counting the representations
of an integer as a sum of $4$ squares, compare \cite[Ch.~11]{Hurwitz}
or \cite[Ch.~XX]{Hardy}, and to observe the relation between
$\lvert\bs{q}\rvert^2$ and $\mathrm{den} (R(\bs{q}))$. This results in
\[
\begin{split}
   f(1) & = 1 \quad \mbox{and} \quad f(2n) = 0 \, , \\
   f(p^r) & = (p\!+\!1)\, p^{r-1} \, , \quad 
          \mbox{for odd primes and $r\geq 1$, and} \\
   f(mn) & =f(m) f(n)\, , \quad 
          \mbox{for $m,n$ coprime (multiplicativity of $f$),}
\end{split}
\]    
see also \cite{Grimmer84,Z2}.  This can be summarized as follows.
\begin{prop}
The Dirichlet series generating function for the number $f(m)$ of 
CSLs of\/ $\ZZ^3$ of index\/ $m$  reads
\begin{equation}  \label{cubgenfun}
\begin{split}
   \varPhi(s) & \; = \; \sum_{m=1}^{\infty} \frac{f(m)}{m^s}
                 \; = \; \prod_{p \neq 2}
                 \frac{1+p_{}^{-s}}{1-p_{}^{1-s}} \; = \;
      \frac{1-2^{1-s}}{1+2^{-s}} \cdot
      \frac{\zeta(s) \zeta(s-1)}{\zeta(2s)}     \\
& \; = \;  1+ \tfrac{4}{3^s} + \tfrac{6}{5^s} + \tfrac{8}{7^s} + 
         \tfrac{12}{9^s}
                + \tfrac{12}{11^s} + \tfrac{14}{13^s} + \tfrac{24}{15^s} 
                + \tfrac{18}{17^s} + \tfrac{20}{19^s} + \tfrac{32}{21^s}
                + \tfrac{24}{23^s} + \tfrac{30}{25^s}+ 
                \cdots  
\end{split}
\end{equation}
\end{prop}

With this generating function, one can again determine the
asymptotic behaviour of $f(m)$. The result is that the
number of CSLs of index $\le N$ is asymptotically given
by $3 N^2/\pi^2$, and the number of
coincidence rotations with index $\le N\/$ by
$72N^2/\pi^2$. 

To expand on the systematics of our generating functions, 
let us remark that Eq.~(\ref{cubgenfun}) can be rewritten as
\be
    \varPhi(s) \; = \; \frac{1}{1+2^{-s}} \cdot
                    \frac{\zeta^{}_H(s/2)}{\zeta(2s)}
\ee
where $\zeta^{}_H(s)=(1-2^{1-2s})\,\zeta(2s)\,\zeta(2s-1)$ is the 
zeta function of the Hurwitz ring $\mathbb{J}$ of integer quaternions, 
compare \cite[$\S$ 63, A.~15 on p.~252]{Scheja}. It is the
generating function for the number of non-zero right ideals
of the ring $\mathbb{J}$, which is a maximal order of the
standard quaternion algebra over $\QQ$.

This is not quite the end of the story.  Coincidence rotations 
of a given index $m\/$ can be collected into equivalence classes
of rotations related by the action of the point symmetry
$\mathcal{O}$ or, more generally, of $\mathrm{O}(3,\RR)$.
This requires a double coset analysis that is described 
in \cite{Grimmer84,Z2}. It turns out that
inequivalent CSLs of $\ZZ^3$  with the same index occur 
for the first time at $\gS=13$.

\begin{figure}[ht] 
\centerline{}
\centerline{\epsfysize=7cm \epsfbox{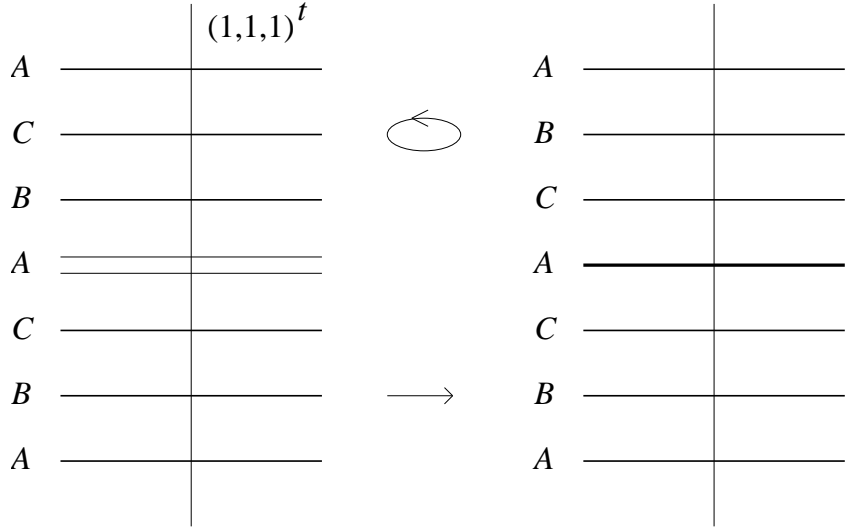}}
\caption{Construction of a twin based on $R(\bs{q})$ with
 $\bs{q} = (0,1,1,1)$. \label{stack}}
\end{figure}

Also, describing the fine structure of a coincidence rotation requires 
an analysis of the lattice planes perpendicular to the rotation axis. 
For example,
the (unique) equivalence class for $\gS=3$ can be represented by
$\bs{q}=(0,1,1,1)$, i.e., a rotation through $\pi$ around $(1,1,1)^t$.
Here, three layers are stacked periodically (with period $\sqrt{3}$), 
with perfect coincidence in one (called $A$, which is therefore an 
ideal grain boundary), but none in the other two 
(called $B,C$, which are interchanged in type by the rotation).
So, this coincidence rotation can be used to create a twin that shows
up as a change in the stacking sequence, see Figure~\ref{stack}.

Let us take a short look at the other cubic lattices, $\gG_F$
and $\gG_B$ (from now on, we use indices $F,P,B$ to distinguish
the three cubic lattices). They possess the same group of 
coincidence isometries as $\gG_P = \ZZ^3$,
\be
      \mathrm{OC}(\gG_F)  \; = \; \mathrm{OC}(\gG_P)
                                \; = \; \mathrm{OC}(\gG_B)
                                \; = \; \mathrm{O}(3,\QQ) \, ,
\ee
but even more, one has
\begin{prop}
   In\/ $\EE^3$, the coincidence index of an orthogonal matrix 
   is the same for all three cubic lattices.
\end{prop}

\noindent {\sc Proof}: 
We can apply Lemma~\ref{sub3}. From the inclusion
\begin{equation}
        \gG_F \; \stackrel{2}{\subset} \;
        \gG_P \; \stackrel{2}{\subset} \;
        \gG_B \, ,
\end{equation}
we know that $\gS_B \, | \, 2 \gS_P$ and $\gS_P \, | \, 
2 \gS_F$. But $\gS_B = \gS_F$ from
Theorem~\ref{sigdual}, because $\gG_B^* = \gG_F^{}$
in this setting. Since $\gS_P$ is always odd, we can either have
$\gS_B=\gS_P$ or $\gS_B=2\gS_P$.

However, as $\gG_B$ and $\gG_P$ have the
same point symmetry, $\gS_P=1$ implies $\gS_B=1$,
and $\gG_B$ cannot possess coincidence
isometries of index 2. From the multiplicativity of the number
of CSLs of index $m$, see \cite{Z2,BLP96}, 
it then follows that no CSL of $\gG_B$
can have even index, hence $\gS_B=\gS_P=\gS_F$. \qed

\smallskip
One can see the second step also without reference to the
multiplicativity of the CSL counting function. Observe that
$\gG_B = \gG_P \,\dot{\cup}\, (\bs{v} + \gG_P)$, with
$\bs{v} = \tfrac{1}{2} (1,1,1)^t$. These two cosets have
disjoint shells, and one can check that coincidence rotations
produce coinciding points in both shells with equal density,
which then implies that $\gS_B = \gS_P$.

Nevertheless, there are specific differences in the coincidence
structure of the three cubic lattices. They show up in different layer
arrangements \cite{Grimmer74}, but details go beyond the scope of
this article.

\subsection{$\bs{d=4}$: the two hypercubic lattices}

In four dimensions, there are only two hypercubic lattices
($P$-type and centred, represented by $\ZZ^4$ and $D_4$),
and they confront us with a different situation that occurs
in no other dimension: the centred lattice, $D_4$, has a larger 
symmetry than its partner of $P$-type, $\ZZ^4$. The dual lattice of 
$D_4$, the weight lattice $D_4^*$, is equivalent to $D_4$ (through 
a rotation followed by a homothety) and need not be considered 
separately, compare \cite{Conway}. 

Let us start with the description of the coincidence structure 
of $D_4$. First, we observe that quaternions give us again a
helpful parametrization of rotations \cite{duVal,Koecher}, where
we now need a pair $(\bs{q}^{}_{1},\bs{q}^{}_{2})$.
The corresponding matrix is defined through 
\[
   M(\bs{q}^{}_1,\bs{q}^{}_2)\bs{x}^t  \; := \;
   \big(\bs{q}^{}_1 \bs{x} \, \overline{\bs{q}}^{}_2 \big)^t  ,
\]   
where $\bs{x}^t$ stands for the transpose of the
quaternion $\bs{x}$ (and is thus a
column vector).  With $\bs{q}^{}_1 = (k,\ell,m,n)$ and $\bs{q}^{}_2 =
(a,b,c,d)$, one finds the $4\!\times\!4$-matrix
\[
M(\bs{q}^{}_{1},\bs{q}^{}_{2}) \; = \;
\mbox{\small $
\left( \!\! \begin{array}{rrrr}
 ak\!+\!b\ell\!+\!cm\!+\!dn  & -a\ell\!+\!bk\!+\!cn\!-\!dm & 
 -am\!-\!bn\!+\!ck\!+\!d\ell & -an\!+\!bm\!-\!c\ell\!+\!dk \\
 a\ell\!-\!bk\!+\!cn\!-\!dm  & ak\!+\!b\ell\!-\!cm\!-\!dn  & 
 -an\!+\!bm\!+\!c\ell\!-\!dk & am\!+\!bn\!+\!ck\!+\!d\ell  \\
 am\!-\!bn\!-\!ck\!+\!d\ell  & an\!+\!bm\!+\!c\ell\!+\!dk  & 
 ak\!-\!b\ell\!+\!cm\!-\!dn  & -a\ell\!-\!bk\!+\!cn\!+\!dm \\
 an\!+\!bm\!-\!c\ell\!-\!dk  & -am\!+\!bn\!-\!ck\!+\!d\ell & 
 a\ell\!+\!bk\!+\!cn\!+\!dm  & ak\!-\!b\ell\!-\!cm\!+\!dn
\end{array} \!\! \right)  $}
\]
which has 
\[
    \det(M(\bs{q}^{}_{1},\bs{q}^{}_{2}) )
    \; = \; (k^2+\ell^2+m^2+n^2)^2\cdot(a^2+b^2+c^2+d^2)^2
           \; = \; \lvert \bs{q}^{}_{1} 
           \rvert^4\,\lvert\bs{q}^{}_{2}\rvert^4
\]     
and  satisfies the orthogonality relation
\[
    M M^t \; = \; \sqrt{\det(M)} \cdot \mathbbm{1}^{}_4 \, .
\]    
Consequently, $|\bs{q}^{}_1| = |\bs{q}^{}_2| = 1$ results in a 4d
rotation matrix, and the group homomorphism $M$: 
$\mathbb{S}^3\times \mathbb{S}^3 \longrightarrow
\mathrm{SO}(4)$ is onto. It provides a twofold cover of the 
rotation group \cite{Koecher},
with $M(\bs{q}^{}_1,\bs{q}^{}_2) = M(-\bs{q}^{}_1,-\bs{q}^{}_2)$.

{}From here, we can find a parametrization of $\mathrm{SO}(4,\QQ)$ if we
start from two non-zero primitive integer quaternions $\bs{q}^{}_1,\bs{q}^{}_2$ 
(i.e., each quaternion has the form $(\kappa,\lambda,\mu,\nu)$
with $\kappa,\lambda,\mu,\nu\in\ZZ$ and $\gcd(\kappa,\lambda,\mu,\nu)=1$).
Then, if we consider the matrix 
\begin{equation}
     R(\bs{q}^{}_1,\bs{q}^{}_2) \; := \;
      M\left(\frac{\bs{q}^{}_1}{|\bs{q}^{}_1|},
         \frac{\bs{q}^{}_2}{|\bs{q}^{}_2|}\right)  \; = \;
         \frac{1}{|\bs{q}^{}_1\bs{q}^{}_2|}
  M(\bs{q}^{}_1,\bs{q}^{}_2) \, ,
\end{equation}
we see that it is an element of $\mathrm{SO}(4,\QQ)$ if and only if 
$|\bs{q}^{}_1\bs{q}^{}_2|^2$
is a {\em square}\/ in $\NN$, in which case we call the pair of integral
quaternions {\em admissible}. But with all admissible pairs of
primitive quaternions, we actually exhaust $\mathrm{SO}(4,\QQ)$, and 
obtain each element twice (due to
$R(\bs{q}^{}_1,\bs{q}^{}_2)=R(-\bs{q}^{}_1,-\bs{q}^{}_2)$).

From Theorem~\ref{cubic1}, we already know that
$\mathrm{(S)OC} (D_4) = \mathrm{(S)O} (4,\QQ)$.  
As in the previous cases, it is sufficient to treat rotations,
since reflections can be written as a product of a rotation
with a special reflection that leaves $D_4$ invariant --
in complete analogy to the situation in the square lattice
($\mathrm{O} (4) = \mathrm{SO} (4)\rtimes C_2$ is a
semi-direct product). So, we need to know the
coincidence index of an arbitrary rotation $R\in\mathrm{SO} (4,\QQ)$.
This is {\em not}\/ just the denominator of $R$, but given
by the following result, see \cite{Z3,BZ} for a proof.

\begin{prop}
Let\/ $(\bs{q}_1^{},\bs{q}_2^{})$ be an admissible pair of
primitive integral quaternions, and let $\gS(\bs{q})$ denote
the index defined above in Eq.~$\eqref{cubind}$.
Then, the matrix $R(\bs{q}_1^{},\bs{q}_2^{}) \in \mathrm{SO}(4,\QQ)$ 
has coincidence index
\begin{equation} \label{d4ind}
  \gS^{}_F(\bs{q}_1^{},\bs{q}_2^{}) \; = \;
      \mathrm{lcm}\left\{ \gS(\bs{q}_1^{}),
                 \gS(\bs{q}_2^{}) \right\} 
\end{equation}
\end{prop}
Here, $\mathrm{lcm}$ denotes the \underline{l}east \underline{c}ommon 
\underline{m}ultiple, and the subscript $F$ refers to the
(face-) centred lattice $D_4$.
With this formula, it is now a combinatorial problem to
determine the number of coincidence rotations of a given
index $m$ and, dividing by 576 (the order of the
rotation symmetry group of $D_4$), also the number
$f^{}_F(m)$ of different CSLs of $D_4$ of index $m$.
As follows explicitly from Eq.~(\ref{d4ind}) and also
from the unique factorization of integral quaternions,
$f^{}_F(m)$ is again a multiplicative function, so that we only
need to calculate it for $m$ a prime power.

Starting from Eq.~\eqref{d4ind} and counting the possibilities to
contribute to $\gS^{}_{F} (p^r)$, it is not difficult to derive
(for $r\ge 1$) the explicit expression
\begin{equation} \label{cubrec}
  f^{}_{F} (p^r) \; = \; f(p^r) \Big( f(p^r) + 
  2 \sum_{\ell=1}^{\left[\tfrac{r}{2}\right]}
                      f(p^{r-2\ell}) \Big)
\end{equation}
with the $f(m)$ of the 3d cubic case in Eq.~(\ref{cubgenfun}). 
An empty sum is to be understood as $0$, and $[.]$ denotes
Gau{\ss}' brackets. The result is (see \cite{Z3} for a proof)
\begin{eqnarray*}
    f^{}_{F} (1) & = &1 \, ,  \\[1mm]
    f^{}_{F} (mn) & = & f^{}_{F} (m) f^{}_{F} (n)\, ,
    \quad \mbox{if $m,n$ coprime}\, , \\[1mm] 
    f^{}_{F} (2m) & = & 0\, , \quad \mbox{and}  \\
    f^{}_{F}(p^r) & = & \frac{p+1}{p-1}\, p^{r-1} 
   \big(p^{r+1} + p^{r-1} - 2 \big) \, ,
    \quad \mbox{for odd primes and $r\geq 1$}\, .
\end{eqnarray*} 
This fixes the Dirichlet series generating function and one obtains
\begin{prop}
   The Dirichlet series generating function for the numbers
   $f^{}_{F} (m)$ of CSLs of index\/ $m$ in the root lattice
   $D^{}_{4}$ reads
\be
\begin{split}
\varPhi^{}_F(s) & \; = \; 
            \sum_{m=1}^{\infty} \frac{f^{}_F(m)}{m^s}
                 \; = \; \prod_{p \neq 2}
                 \frac{(1+p_{}^{-s}) (1+p_{}^{1-s})}
                      {(1-p_{}^{1-s})(1-p_{}^{2-s})}  \\
& \; = \; 1+\tfrac{16}{3^s}+\tfrac{36}{5^s}+\tfrac{64}{7^s}+\tfrac{168}{9^s}+
\tfrac{144}{11^s}+\tfrac{196}{13^s}+\tfrac{576}{15^s}+\tfrac{324}{17^s}+
\tfrac{400}{19^s}+\tfrac{1024}{21^s}+\tfrac{576}{23^s}+
\cdots   
\end{split}
\ee
\end{prop}
Again, the possible indices are precisely all odd integers,
$\gS_{D_4} (\mathrm{SO} (4,\QQ)) = 2\NN_0 + 1$.
A comparison with the cubic case in $3$-space reveals the remarkable
identity
\begin{equation} \label{recursion}
      \varPhi^{}_F (s) \; = \; \varPhi (s) \, \varPhi (s-1) 
\end{equation}
where $\varPhi(s)$ is the generating function of  Eq.~(\ref{cubgenfun}). 
This makes is rather easy to calculate the asymptotic behaviour
from that in $3$ dimensions. The right-most pole of $\varPhi^{}_{F} (s)$
is at $s=3$, so one obtains that the number of CSLs of index $\le N$
grows asymptotically as $\frac{210}{\pi^6}\zeta(3) N^3 \simeq 0.26257
\, N^3$ (note that $\zeta(3)$ is known to be irrational, but its value is
only known numerically). For the number of coincidence rotations
with index $\le N$, one has to multiply by $576$.
Another consequence of Eq.~\eqref{recursion} is the formula
\begin{equation}
   f^{}_F(m) \; = \; \sum_{d\,|\,m} d\cdot f(d)
                     \cdot f(m/d)
\end{equation}
which follows from the convolution theorem of Dirichlet series 
and allows for a simple and efficient calculation of the numbers
$f^{}_{F} (m)$.

Having described the root lattice $D_4$, the (face-)centred cubic lattice
in four dimensions, we now turn to the slightly more complicated 
case of the $P$-type cubic lattice, $\ZZ^4$. From the inclusion
\be
     D_4^{}   \; \stackrel{2}{\subset} \;   \ZZ^4 
                  \; \stackrel{2}{\subset} \;   D_4^*
\ee
and the result that $D_4^{}$ and $D_4^*$ have the same
index formula, see Theorem~\ref{sigdual},
we get $\gS_F \, | \, 2\gS_P$ and
$\gS_P \, | \, 2\gS_F$. But $\gS_F$ is always odd,
so either $\gS_P = \gS_F$ or $\gS_P = 2\gS_F$.
Here, in contrast to the situation in $3$-space, both possibilities arise.
The point symmetry group of $D_4$ is larger than that of $\ZZ^4$, with
\be
      [ \mbox{Aut}(D_4) : \mbox{Aut}(\ZZ^4) ] \; = \; 3   \,.
\ee
Consequently, one third of the elements of $\mbox{Aut}(D_4)$
are symmetries of $\ZZ^4$ while the others result in
coincidence isometries of $\ZZ^4$. They turn out to have index $2$. 
In going from here to the number of CSLs of index $m$, it
is clear that we have $f_P^{}(1)=1$,
$f_P^{}(2) = 2$ and $f_P^{}(2^r)=0$ for $r>1$.
The multiplicativity (which needs to be proved, e.g., 
similarly to the arguments given in \cite{Z3})
of $f_P^{}(m)$ then gives the general answer
\[
\begin{split}
  f_P^{}(m) & \; = \; f^{}_F(m)\, , \qquad \;\, \mbox{for $m$ odd,} \\
  f_P^{}(m) & \; = \; 2f^{}_F(m/2)\, , \quad 
  \mbox{for $m \equiv 2$ (4), and} \\
  f_P^{}(4m) & \; = \; 0 \, .
\end{split}
\]
This can now easily be summarized as follows.
\begin{prop}
The Dirichlet series generating function for the number\/ $f_P (m)$ 
of CSLs of index\/ $m$ in\/ $\ZZ^4$ reads
\be
\begin{split}
    \varPhi^{}_P(s) & \; = \; (1+2^{1-s})\cdot\varPhi_F(s) \; = \;
       (1+2^{1-s}) \prod_{p \neq 2}
                 \frac{(1+p_{}^{-s}) (1+p_{}^{1-s})}
                      {(1-p_{}^{1-s})(1-p_{}^{2-s})}  \\
& \; = \; 1+\tfrac{2}{2^s}+\tfrac{16}{3^s}+\tfrac{36}{5^s}+\tfrac{32}{6^s}+
\tfrac{64}{7^s}+\tfrac{168}{9^s}+\tfrac{72}{10^s}+\tfrac{144}{11^s}+
\tfrac{196}{13^s}+\tfrac{128}{14^s}+\tfrac{576}{15^s}+\tfrac{324}{17^s}
+\cdots
\end{split}
\ee
\end{prop}

Note that, in comparison to the $D_4$ case, the number of CSLs grows
faster by a factor of $5/4$, while the number of coincidence rotations
grows slower by a factor of $5/12$, due to the smaller point symmetry
group of $\ZZ^4$.

This calculation was actually possible without giving the
corresponding index formula first, by using the multiplicativity of
$f^{}_P(m)$. Since $\mathrm{SOC}(\ZZ^4)=
\mathrm{SOC}(D_4)=\mathrm{SO}(4,\QQ)$, there must be a different index
formula for $\ZZ^4$, and indeed one obtains, for an admissible pair of
primitive integral quaternions (see \cite{Z3} for a proof)
\begin{equation}
\begin{split}
   \gS^{}_P(\bs{q}_1^{},\bs{q}_2^{}) = &
      \mbox{lcm}\left\{ \gS(\bs{q}_1^{}),
                 \gS(\bs{q}_2^{}),
      \mbox{den}(R(\bs{q}_1^{},\bs{q}_2^{})) \right\} \\
  = & \mbox{lcm}\left\{\gS^{}_F(\bs{q}_1^{},\bs{q}_2^{}),
       \mbox{den}(R(\bs{q}_1^{},\bs{q}_2^{})) \right\} \; .
\end{split}
\end{equation}
At this point, we close our description of the coincidence structure of
lattices and turn to the perhaps more interesting case of quasicrystals.

\section{Coincidence isometries for modules}
 
So far, we have described the case of lattices and crystals. For the 
treatment of quasicrystals, we have to extend our concepts to 
$\ZZ$-modules, embedded in Euclidean space, 
which are not necessarily discrete any more.
\begin{defin}  \label{mod1}
  A subset $\mathcal{M}$ of $\EE{\hspace{0.2pt}}^d$ is called a 
  {\em $\ZZ$-module}, of rank $r$ and dimension $d$,
  when it is the $\ZZ$-span of $r$ vectors
  $\bs{a}_1^{}, \ldots, \bs{a}_r$ (the {\em basis} of the module)
  that are linearly independent over $\ZZ$, 
  but span $\EE{\/}^d$ over $\RR$.
\end{defin}
Clearly, we must have $r\geq d$, and a module with $r=d$ is a lattice.
As a group, a module of rank $r$ is isomorphic to the free Abelian group
of rank $r$ -- so we may consider such modules as special geometric
realizations of free Abelian groups.

At this point, the concepts of {\em submodule}, {\em
  commensurateness}, and that of {\em coincidence isometry\/} and {\em
  index}, are defined in exact analogy to Section~\ref{sec2}, so there
is no need to repeat them here.  We only note that the set of
coinciding points is now a module rather than a lattice, wherefore we
call it a {\em \underline{c}oincidence \underline{s}ite
  \underline{m}odule}, or CSM for short.  A bit more care is needed
for the definition of a {\em dual\/} module.  Let us first observe
\begin{lemma}
  {}For every module $\mathcal{M} \subset \EE{}^d$ of rank $r$ and 
  dimension $d\leq r$, there is a lattice $\gG\subset\EE{}^r$ 
  such that $\mathcal{M}$
  is the one-to-one projection of $\gG$ into $\EE{}^d$.
\end{lemma}

\noindent {\sc Proof}: 
Out of the basis $\{\bs{a}_1^{}, \ldots, \bs{a}_r\}$ of $\mathcal{M}$, 
we can pick $n$ $\RR$-linearly independent vectors which span
$\EE{}^d$ over $\RR$, say $\bs{a}_1^{}, \ldots, \bs{a}_d$ w.l.o.g.
They span a lattice in $\EE{}^d$. The statement is now obvious as
we can add one new dimension for each basis vector remaining.  
This can always be done in such a way that the projection is
one-to-one on $\gG$.     \qed

\smallskip
It is an obvious idea to try to define a dual object for a module 
through a lift to a lattice $\gG$ because then $\gG^*$ is
well-defined and can be projected down again. Unfortunately, this
is neither unique nor satisfactory, as it can happen that the
object defined this way is a module of smaller rank than the 
original one. If, however, there is some additional structure
(e.g., irreducible symmetry), such a lift can be made essentially
unique and the dual object is well-defined, see \cite{Patera,PABP}
for relevant examples in our present context. If this situation
applies, it is again true that a module and its dual module
share coincidence group and index formula.

In what follows, we concentrate on examples that are connected
with the golden ratio, $\tau=(1+\sqrt{5}\,)/2$. The modules that
will appear are invariant under multiplication by $\tau$ or a power
thereof. So, it might be instructive to see how the construction
of a dual module works here. The simplest example is the ring
of golden integers
\be
    \ZZ[\tau] \; := \; \{ m+n\tau \mid m,n\in\ZZ\}\, ,
\ee
which is the ring of algebraic integers in the quadratic field
$\QQ(\tau)$, but can also be seen as a $\ZZ$-module $\mathcal{M}$ of
rank $2$ and dimension $1$ in $\EE$. The special structure that helps
in this case is the existence of an automorphism $'$ of $\QQ(\tau)$,
called {\em algebraic conjugation}, which maps $\tau$ to its algebraic
conjugate $\tau'=-1/\tau = 1-\tau$ (hence the name) and thus $\alpha=a+b\tau$
to $\alpha'=a+b\tau'=(a+b)-b\tau$.  Now, the set \be \gG \; := \; \{
(\alpha,\alpha') \mid \alpha\in\ZZ[\tau]\} \ee is a lattice in
$2$-space $\EE^2$, from which one obtains $\ZZ[\tau]$ by projection
into $\EE$. In $\EE^2$, $\gG$ has a well-defined dual lattice,
$\gG^*$, and neither $\gG$ nor $\gG^*$ has a lattice direction
parallel to $\EE$. But then, we {\em define\/} the dual module
$\mathcal{M}^*$ by the projection of $\gG^*$ into $\EE$
(which, in this case, gives $\mathcal{M}^*=\mathcal{M}/\sqrt{5}
=\ZZ[\tau]/\sqrt{5}\,$), and this object is unique in the sense
that any other embedding of $\mathcal{M}$ into $\EE^2$ which
maps algebraic conjugation to a lattice automorphism will result
in the same dual object $\mathcal{M}^*$.

Though examples in higher dimensions are more complicated (and
will require more than just one automorphism), the basic idea is
similar. In particular, it applies to all quasicrystals of
interest, compare \cite{PABP}. It turns out that root lattices
\cite{Conway} prove extremely handy here \cite{BJKS}.
Let us now illustrate the coincidence problem for
noncrystallographic patterns by a seris of examples
related to fivefold symmetry. For further material on 
planar structures with $N$-fold symmetry, we refer to
\cite{BP94,Pleasants,BG}.

\section{Modules and quasicrystals: the root systems 
         $H_2$, $H_3$ and $H_4$}

\begin{figure}[ht] 
\centerline{}
\centerline{\epsfysize=15mm
\epsfbox{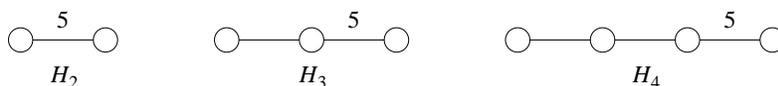}}
\caption{The non-crystallographic Coxeter groups of type $H$. 
  \label{hroot}}
\end{figure}

Among the many possible quasiperiodic tilings, those attached
to fivefold symmetry are of particular interest, especially in
view of their application in solid state physics. Let us therefore
start from the exceptional Coxeter groups $H_2$ (usually called
$I_2(5)$), $H_3$ and $H_4$ \cite{Humphreys} shown in 
Fig.~\ref{hroot}. They are the symmetry groups of certain
regular polytopes \cite{Coxeter}, namely of the regular
decagon ($\{10\}$ in Schl\"afli's notation), 
the icosahedron $\{3,5\}$ (or dodecahedron $\{5,3\}$) and the regular
600-cell $\{3,3,5\}$ (or the regular 120-cell $\{5,3,3\}$),  
and are of order 20, 120 and 14400, respectively.

The corresponding root systems (up to normalization) are given by the
vectors which point to the 10 vertices of a regular decagon ($H_2$),
to the 30 vertices of the icosidodecahedron 
$\big\{\begin{smallmatrix}  3 \\ 5
\end{smallmatrix}\big\}$ ($H_3$), and to the 120 vertices of the regular
600-cell in 4-space ($H_4$), see \cite{Patera} or \cite{Humphreys} for
details. The $\ZZ$-spans of these root systems define $H_n$-symmetric
modules, which can be seen as projections from the root lattices
$A_4$, $D_6$ and $E_8$. In particular, they all have well-defined
duals, though we do not expand on this question here -- it is
discussed in detail in \cite{Patera}.

Having set the scene, we can now describe the coincidence structure of these
modules (and some closely related ones). We shall also briefly discuss the
connection to quasiperiodic tilings.

\subsection{$\bs{d=2}$: Coincidence rotations for tenfold symmetry}

\begin{figure}[ht] 
\centerline{\epsfysize=100mm
\epsfbox{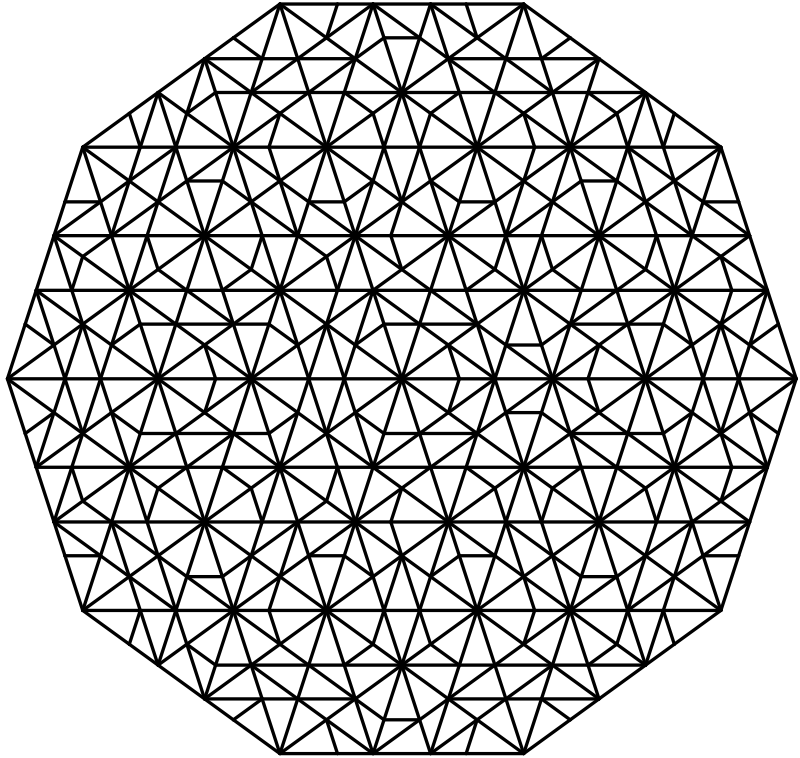}}
\caption{Cartwheel version of the T\"ubingen 
            triangle tiling. \label{trifig}}
\end{figure}

Let us consider a 2d quasicrystal with tenfold symmetry, the 
T\"ubingen triangle tiling \cite{BKSZ} of Fig.~\ref{trifig}, say. 
As mentioned earlier, the coincidence problem
splits into two parts: first, the coincidence problem for the underlying
$\ZZ$-module $\mathcal{M}_{10}$ (which is the limit translation module
\cite{BS} of the tiling) and second, the correction, due to the 
acceptance domain, of the coincidence indices obtained in this way 
\cite{Pleasants}. Here, we discuss in detail only the first part.

In the complex plane, the tenfold module of rank 4 over $\ZZ$ \cite{Mermin}
can be written as the direct sum \be \mathcal{M}_{10} \; = \; \ZZ \, 1 \,
\oplus \, \ZZ \, \xi \, \oplus \, \ZZ \, \xi^2 \, \oplus \, \ZZ \, \xi^3 , \ee
with $\xi = e^{2\pi i/5}$ a {\em fifth} root of $1$.  As such,
$\mathcal{M}_{10}$ is the ring of algebraic integers in the cyclotomic field
$K_5=\QQ(\xi)$ (a field extension of $\QQ$ of degree $4$).  Again, prime
factorization is unique up to units \cite{Wash}. One can now go through 
essentially the same argument as in the case of the square lattice.
The primes in $\ZZ[\xi]$ are slightly more complicated, but only those
dividing a rational prime $p\equiv 1$ $\, (5)$ can enter the
factorization of $e^{i\varphi}=\alpha/\beta$ with $\alpha,\beta
\in\ZZ[\xi]$ coprime. The result is \cite{Pleasants}:
\begin{prop}
Every coincidence rotation $($written as $e^{i\varphi}\in K_5 =
\QQ(\xi))$ of the standard tenfold symmetric module $\mathcal{M}_{10}$
can be factorized as
\begin{equation} \label{factor3}
   e^{i\varphi} \, = \,  \, \varepsilon \cdot \!\!\!
                \prod_{p \equiv 1 \; (5)}
\left( \frac{\omega_p^{(1)}}{\overline{\omega}_p^{(1)}} 
       \right)_{ }^{n_p^{(1)}}
\left( \frac{\omega_p^{(2)}}{\overline{\omega}_p^{(2)}} 
       \right)_{ }^{n_p^{(2)}}
\end{equation}
where $\varepsilon$ is a\/ $10$th root of\/ $1$ and thus a
unit in\/ $\ZZ[\xi]$, and only finitely many of the exponents
$n^{(1)}_{p}$, $n^{(2)}_{p}$ are different from\/ $0$.
\end{prop}
 
This factorization is slightly more complicated than that in the case of
the square lattice, as we have two {\em independent\/} 
generators for each basic index $p \equiv 1 \; (5)$. This originates
in the fact that these (rational) primes are the product of $4$
cyclotomic primes in $\ZZ[\xi]$ which form two independent pairs
of complex conjugates.

The index of a rotation $R$ is
\begin{equation}
    \gS(R) \, = \, \prod_{p \equiv 1 \; (5)}
                             p_{ }^{(|n_p^{(1)}| + |n_p^{(2)}|)} 
\end{equation}
and the group of coincidence rotations has the form
\begin{equation}
 \mathrm{SOC}(\mathcal{M}_{10}) \; \simeq \; C_{10} 
  \times \ZZ^{(\aleph_0)} \, .
\end{equation}
Thus, in spite of the more complicated factorization,
the structure of the coincidence group remains simple
(as it does for other planar symmetries, compare \cite{Pleasants,BG}).
Let us give two examples before we continue. 
$\gS=11$ and $\gS=31$ are the smallest
non-trivial indices, with two generators each. 
With $\xi=e^{2\pi i/5}$, they can be written as follows:
\be
   \gS=11 \,:  \;\;\; \frac{2+\xi}{2+\overline{\xi}} \; , \;\;\;  
             \frac{2+\xi^2}{2+\overline{\xi}^2} \; ; \hspace*{8mm}
   \gS=31 \,:  \;\;\; \frac{2-\xi}{2-\overline{\xi}} \; , \;\;\;  
             \frac{2-\xi^2}{2-\overline{\xi}^2} \; .
\ee
Coincidence reflections can be described as products of rotations
with the reflection in the $x$-axis, exactly as in the case of the
square lattice, so we need not repeat that argument here.

If $10 f(m)$ denotes the number of coincidence rotations of index $m$,
the multiplicative function $f(m)$ actually counts the number of 
different CSMs of index $m$. We then obtain, with
$\QQ (\tau) = \QQ(\sqrt{5}\,) = \QQ(\xi)\cap \RR$:

\begin{prop}
The Dirichlet series generating function for the coincidence problem
of the tenfold symmetric module of rank\/ $4$ in the plane reads
\be
\begin{split}
 \varPhi(s) & \; = \; \sum_{m=1}^{\infty} \frac{f(m)}{m^s}
                 \; = \, \prod_{p \equiv 1 \; (5)}
           \left( \frac{1+p_{}^{-s}}{1-p_{}^{-s}} \right)^2     
              \;=\;  \frac{1}{1+5^{-s}}\cdot
          \frac{\zeta^{}_{\QQ(\xi)}(s)}{\zeta^{}_{\QQ(\tau)}(2s)}  \\
& \; = \; 1+\tfrac{4}{11^s}+\tfrac{4}{31^s}+\tfrac{4}{41^s}+\tfrac{4}{61^s}+
\tfrac{4}{71^s}+\tfrac{4}{101^s}+\tfrac{8}{121^s}+\tfrac{4}{131^s}+
\tfrac{4}{151^s}+\tfrac{4}{181^s} + \cdots
\end{split}
\ee
\end{prop}
Here, $\zeta^{}_{\QQ (\xi)} (s)$ is the Dedekind zeta function \cite{Wash} 
of the cyclotomic field $\QQ (\xi)$, while $\zeta^{}_{\QQ(\tau)}$ is the
zeta function of the maximal real subfield, $\QQ (\xi + \bar{\xi})=
\QQ(\tau)$, see also Eq.~\eqref{xxx} below.

All CSMs are scaled versions of $\mathcal{M}_{10}$ 
and the number of CSMs of index $\le N\/$ is asymptotically
$5 \log(\tau) N/\pi^2$ (while the number of coincidence rotations with
index $\le N$ is $10$ times as large).

How do these results apply to the coincidence problem of the
tenfold symmetric triangle tiling? The latter can be
obtained through projection from the root lattice $A_4$ 
with a regular decagon as its window \cite{BKSZ}. 
A coincidence in the set of vertex points occurs if and only if
there is a coincidence in the module $\mathcal{M}_{10}$ such that the
image point in internal space lies both in the original window
{\em and\/} in an appropriately rotated window.
A consequence of this is that the coincidence group of the tiling is 
still $\mathrm{SOC}(\mathcal{M}_{10})$, but also that the $\gS$-factor 
or ``index'' of  each group element is normally smaller
than its index in $\mathcal{M}_{10}$ by a correction factor close to $1$ 
(depending on the group element).  

This is called the {\em window correction factor}, and has to be calculated
for each tiling separately, as a function of the rotation angle.
It also explains why the set of coinciding points forms
a tiling of slightly different type from the original one, a small 
proportion of the points of the original tiling being missing from it.  
In fact the term ``index" for the reciprocal of the fraction of coinciding 
points is no longer
appropriate in this setting, as it neither has a purely algebraic
interpretation, nor is an integer any more.  Details of this and 
the determination of the rotation angle in internal space by means of 
algebraic conjugation are given in \cite{Pleasants}.
In many examples, in particular when the windows are regular
polytopes, the maximal error is so small (at most of the order
of a few percent) that one can safely ignore it and work with
the module index instead.

\subsection{$\bs{d=3}$: the icosahedral modules of rank 6}

Icosahedral quasicrystals are of particular interest, and one would like to
know their coincidence structure in detail \cite{WaLu,Radu95}.  We restrict
our discussion to the investigation of the 3 different 3d icosahedral modules
of rank 6 over $\ZZ$ \cite{Rokhsar} and again omit the determination of the
window correction.  We shall call the modules $\mathcal{M}_B$,
$\mathcal{M}_P$, $\mathcal{M}_F$ for $B$-, $P$- and $F$-type,
respectively\footnote{This terminology originates from the fact that these
  modules can be obtained as projections of the three types of hypercubic
  lattices in $6$-space, $D^{*}_{6}$, $\ZZ^6$, and $D^{}_{6}$.}. They are
spanned by the orthonormal basis $\bs{e}_1,\bs{e}_2,\bs{e}_3$ with
coefficients $\alpha_i \in \ZZ[\tau]$, $\tau = (1 + \mbox{\small
  $\sqrt{5}\,$})/2$, as follows:
\begin{eqnarray}
   \mathcal{M}_B & = & \{\, \mbox{\small $\sum$}_{i=1}^3 \; 
                               \alpha_i \bs{e}_i \; | \;
       \tau^2 \alpha_1 + \tau \alpha_2 + \alpha_3 \equiv 0 \;\; (2) \, 
                               \}  \nonumber \\
                               \label{icomod}
   \mathcal{M}_P & = & \{\;\; \bs{x} \in \mathcal{M}_B   \;\, | \;
       \alpha_1 + \alpha_2 + \alpha_3 \equiv 0 \; \mbox{ or } \; 
                               \tau \;\; (2) \, \} \\ 
   \mathcal{M}_F & = & \{\;\; \bs{x} \in \mathcal{M}_B   \;\, | \; 
        \alpha_1 + \alpha_2 + \alpha_3 \equiv 0 \;\; (2) \, \} \; .  
                               \nonumber
\end{eqnarray}
The use of an orthonormal basis may be a bit surprising at first sight, but it
will prove useful in a moment. It is possible in this simple form because
$\bs{e}^{}_{1},\bs{e}^{}_{2},\bs{e}^{}_{3}$ are chosen parallel to $3$
mutually orthogonal twofold axes of the icosahedron.  In this setting,
$\mathcal{M}_F$ is the $\ZZ$-span of the root system of type $H_3$, and hence
a $\ZZ$-module of rank $6$ and dimension $3$. It actually also is a
$\ZZ[\tau]$-module (of rank $3$), which is also true of $\mathcal{M}_B$, but
not of $\mathcal{M}_P$. This relates to the fact that $\mathcal{M}_F$ and
$\mathcal{M}_B$ are invariant under multiplication by $\tau$, while
$\mathcal{M}_P$ is only invariant under multiplication by $\tau^3$.

To describe the coincidence rotations (in the orthogonal basis),
Cayley's parametrization can again be used.  Our first
assertion is that the coincidence group is the same for all three
modules:
\begin{prop} $\;
 \mathrm{OC}(\mathcal{M}_B) \,=\, \mathrm{OC}(\mathcal{M}_P) \,=\, 
 \mathrm{OC}(\mathcal{M}_F) \,=\,
 \mathrm{O}(3,\QQ(\tau)) \; .$
\end{prop}

\noindent {\sc Proof}: 
Observe the relations
$2\ZZ[\tau]^3 \stackrel{4}{\subset} \mathcal{M}_F \stackrel{2}{\subset}
\mathcal{M}_P \stackrel{2}{\subset} \mathcal{M}_F \stackrel{4}{\subset} 
\ZZ[\tau]^3$ and $\ZZ[\tau]^3=\ZZ^3 \oplus\tau\ZZ^3$. 
These modules all possess the same $\mathrm{OC}$-group, and this obviously is 
$\mathrm{O}(3,\QQ(\tau))$ by the same type of argument we have used 
previously in the
discussion of the (hyper-)cubic lattices.   \qed

\smallskip
The unit quaternions $(1,0,0,0)$, $\mbox{\footnotesize 
$\tfrac{1}{2}$}(1,1,1,1)$,
$\mbox{\footnotesize $\tfrac{1}{2}$} (\tau,1,-1/\tau,0)$ together
with all even permutations and arbitrary sign flips form a group
$\mbox{\small $\widehat{Y}$}$ of order 120 
which is the usual double cover \cite{Humphreys,Moody} 
of the icosahedral group 
$Y=\{R\in \mathrm{SO}(3,\QQ(\tau)) \,|\, \gS(R)=1\}$.
The icosian ring $\II$, see \cite{Moody} for details, 
consists of all integral linear combinations of elements
in $\mbox{\small $\widehat{Y}$}$ and is a maximal order with unique 
(left- or right-) factorization in the quaternion algebra over
the field $\QQ(\tau)$. One finds the relation
$\mathrm{SO}(3,\QQ(\tau)) = \{R(\bs{q}) \,|\, \bs{0}\neq\bs{q}\in\II\,\}$, 
and our second assertion is the index formula for a coincidence 
rotation $R_0 \in \mathrm{SO}(3,\QQ(\tau))$, again for all three modules:
\begin{equation}
   \gS(R_0) \; = \; \gcd \, 
    \{\, N( | \bs{q} |^2 ) \; \; | \; \; \bs{q}\in\II \, ,\;  
         R(\bs{q})=R_0  \}\, ,
\end{equation}
where the argument $| \bs{q} |^2$ on the right hand side is 
always a number 
in $\ZZ[\tau]$ and its norm is defined by $N(m+n\tau)=m^2+mn-n^2$. 
We use the convention that the gcd is always a positive number.
The indices $\gS$ run through all positive integers of the form 
$m^2+mn-n^2$ with integral $m\/$ and $n$.
These are the numbers all of whose prime factors congruent to 
2 or 3 (mod 5)
occur with even exponent only. (They can also be characterized as the
positive numbers of the form $5x^2-y^2$ with integral $x\/$ and $y$,
as used in \cite{Radu95}.)
For $\gS \leq 100$, one finds the list of numbers 
\[ 1,4,5,9,11,16,19,20,25,29,31,36,41,44,45,49,
55,59,61,64,71,76,79,80, 81,89,
95,99,100,
\]
which covers the cases known from \cite{WaLu}.

If $60f(m)$ is the number of coincidence rotations of index $m$, 
$f(m)$ is the number of different CSMs of index $m$. Since the
icosian ring is a maximal order with unique (left- or right-)
factorization \cite{Reiner,Vigneras}, $f(m)$ is again
a multiplicative function, i.e., $f(1)=1$ and $f(mn)=f(m)f(n)$
whenever $m,n$ are coprime. Furthermore, with $r\geq 1$, one finds 
\[
     f(5^r) \; = \; 6\cdot 5^{r-1} .
\]     
Then, if $p\equiv\pm 2$ (5), 
\[
   f(p^{2r-1}) \; = \; 0 \quad \mbox{and} \quad
   f(p^{2r}) \; = \; (p^2+1) p^{2(r-1)}.
\]
Finally, if $p\equiv\pm 1$ (5),
\[
    f(p^r) \; = \; (p+1) \big((r+1)p^{r-1} + (r-1) p^{r-2}\big).
\]
This fully determines the generating function 
of $f(m)$.
\begin{prop}
The Dirichlet series generating function for the number of CSMs
of an icosahedral module from Eq.~$\eqref{icomod}$ is given by
\[
\begin{split}
  \varPhi(s) &  \; = \; \sum_{m=1}^{\infty} \frac{f(m)}{m^s}
                 \; = \; \frac{1+5_{}^{-s}}{1-5_{}^{1-s}}
                         \prod_{p \equiv \pm 2 \; (5)}
                 \frac{1+p_{}^{-2s}}{1-p_{}^{2(1-s)}}
                         \prod_{p \equiv \pm 1 \; (5)}
           \left( \frac{1+p_{}^{-s}}{1-p_{}^{1-s}} \right)^2     \\
& \; = \; 1+\tfrac{5}{4^s}+\tfrac{6}{5^s}+\tfrac{10}{9^s}+
\tfrac{24}{11^s}+\tfrac{20}{16^s}
 +\tfrac{40}{19^s}+\tfrac{30}{20^s}+\tfrac{30}{25^s}+\tfrac{60}{29^s}
 +\tfrac{64}{31^s}+\tfrac{50}{36^s}+ \cdots
\end{split}
\]
\end{prop}
The number of CSMs of index $\le N\/$ is 
asymptotically $45\sqrt{5} \log(\tau)N^2/ 2 \pi^4$ (while the number
of coincidence rotations with index $\le N$ is $60$ times as large).

The function $\varPhi(s)$ can be expressed in terms of zeta functions as
\begin{equation} 
 \varPhi(s) \; = \; \frac{\zeta^{}_L(s) \zeta^{}_L(s-1)}{\zeta^{}_L(2s)}
         \; = \; \frac{\zeta^{}_{\II}(s/2)}{\zeta^{}_L(2s)}
\end{equation}
with the quadratic field $L:=\QQ(\tau)=\QQ(\sqrt{5}\,)$, $\zeta^{}_{\II}(s) = 
\zeta^{}_L(2s) \zeta^{}_L(2s-1)$ the $\zeta$-function
of the icosian ring, and
\begin{equation} \label{xxx}
  \zeta^{}_L(s) \; = \; \frac{1}{1-5^{-s}}
             \prod_{p \equiv \pm 2 \; (5)}
             \frac{1}{1-p^{-2s}}
             \prod_{p \equiv \pm 1 \; (5)}
             \frac{1}{(1-p^{-s})^2} ,
\end{equation}
A more complete description of the icosahedral case, which requires
a certain level of mathematical machinery and various results from
algebraic number theory, will be given in \cite{BLP96}.

\subsection{Short digression on a related cubic module of rank 6}

Crystals and quasicrystals are specific examples of ordered phases
with longe-range order, but there are many other ones. Incommensurate
structures are also widely studied in the literature, and they are
also connected to modules rather than lattices. It is an obvious
question what happens to such modules if they still have {\em cubic}
symmetry (and hence dimension 3), but rank 6. There are 6 types of
such modules, called $B+B$, $P+P$, $F+F$, $B+P$, $B+F$, and $P+F$
in a suggestive notation. Let us consider the case $P+P$ in a bit
more detail. It is clear that $\ZZ^3 + \alpha \ZZ^3$ is an example
of it if we demand $\alpha\not\in\QQ$ 
(otherwise, the rank would not be 6).

So, let us assume that $\alpha$ is irrational. Then, $\ZZ^3\cap
\alpha\ZZ^3 = \{0\}$ and we can write the module as
$\ZZ^3 \oplus \alpha \ZZ^3$.
Generically, the $\mathrm{OC}$-group will be that of $\ZZ^3$ itself (with the
index squared), because, if $\alpha$ is not algebraic, there is no
rotation which brings a point of $\ZZ^3$ into coincidence with
one of $\alpha\ZZ^3$. This changes quite a bit if $\alpha$ is
algebraic. A case of particular interest in our context is that
of $\alpha=\tau$, where we get the module
\[
 \mathcal{M}_C \; = \; \ZZ^3 \oplus \tau\ZZ^3 \; = \; \ZZ[\tau]^3 \; .
\]
This module has the same $\mathrm{OC}$-group as the three icosahedral modules
above, $\mathrm{OC}(\mathcal{M}_C)=\mathrm{O}(3,\QQ(\tau))$, 
but a different index formula.
In fact, in complete analogy to the cubic lattices, one can prove that
\be \label{modcubind}
  \gS(R) \; = \; | N( \den(R) ) |  \; .
\ee
Here, $\den(R)$ is defined w.r.t.\ $\ZZ[\tau]$ and hence a number of
the form $m+n\tau$, $N$ is the norm in $\ZZ[\tau]$ as used above,
and the absolute value is needed because the denominator is
only defined up to units in $\ZZ[\tau]$, which are the numbers
$\pm\tau^r$ with norm $(-1)^r$.

It is clear from Eq.~(\ref{modcubind}) that the set of indices is
the same as in the icosahedral case, i.e., all positive integers
which are representable by the integral quadratic form
$m^2 + mn - n^2$. More surprising is the result that the
generating function is very similar: differences only occur
for indices which are divisible by 4. Explicitly:
\be
\begin{split}
 \varPhi_C(s) & \; = \;  
    \frac{1+4^{1-s}}{1+4^{-s}} \cdot
    \frac{\zeta^{}_L(s) \zeta^{}_L(s-1)}
         {\zeta^{}_L(2s)}   \\
& \; = \; 1+\tfrac{8}{4^s}+\tfrac{6}{5^s}+\tfrac{10}{9^s}+
\tfrac{24}{11^s}+\tfrac{32}{16^s}
 +\tfrac{40}{19^s}+\tfrac{48}{20^s}+\tfrac{30}{25^s}+\tfrac{60}{29^s}
 +\tfrac{64}{31^s}+\tfrac{80}{36^s}+ \cdots
\end{split}
\ee

Due to the prefactor, the number of CSMs of index $\le N$ grows faster
than that of the icosahedral case by a factor of $10/9$, while, due  to the
different point symmetry groups, the number of coincidence rotations
grows slower by a factor of $4/9$.

\subsection{$\bs{d=4}$: the icosian ring as $H_4$-symmetric module}

Our final example is the icosian ring itself, compare \cite{Moody},
viewed as the module obtained as the $\ZZ$-span of the root system
of type $H_4$. It is the limit translation module of a highly symmetric
4d quasicrystal studied in \cite{Elser}, and can also be obtained
by projection of the root lattice $E_8$ to a 4d subspace that is
invariant under the action of the symmetry group of the regular
600-cell (which, in turn, is isomorphic to the Coxeter group $H_4$),
see \cite{BM,Moody94} for more.

The solution of this case proceeds in close analogy to that of $D_4$,
details will appear in \cite{BZ}. The coincidence group 
obviously is $\mathrm{OC}(\II) = \mathrm{O}(4,\QQ(\tau))$, and in
order to obtain a parametrization of its rotation matrices, we are
again using admissible pairs of quaternions $\bs{q}^{}_1,\bs{q}^{}_2
\in\II$. This now means that $|\bs{q}^{}_1 \bs{q}^{}_2|^2$ must be
a {\em square} in $\ZZ[\tau]$. This gives some slight extra complication
from the arithmetic of $\ZZ[\tau]$, compare \cite{Dodd}, as
rational primes $p\equiv\pm 1$ (5) split into a product of two
$\ZZ[\tau]$-primes which are algebraic conjugates of one another,
but not associates. Consequently, Eq.~(\ref{cubrec}) has a counterpart
in the present case that has to be modified for such primes
by an extra factor 1/2 on the right hand side, while it is true
in its unaltered form for all other primes. The result is a
multiplicative function $f^{}_{\II}(m)$ which counts the CSMs of
$\II$ of index $m$. It is specified by $f^{}_{\II}(1)=1$, and,
for $r\geq 1$, by
\[
   f^{}_{\II}(5^r) \; = \; \frac{3}{2} 5^{r-1} (5^{r+1} + 5^{r-1} - 2);
\]    
if $p\equiv\pm2$ (5), one has 
\[
    f^{}_{\II} (p^{2r-1}) \; = \; 0  \quad \mbox{and} \quad
    f^{}_{\II} (p^{2r}) \; = \; \frac{p^2+1}{p^2-1} p^{2(r-1)}\,
      \big(p^{2(r+1)} + p^{2(r-1)} - 2 \big).
\]      
Finally, for $p\equiv\pm1$ (5), one obtains the rather lengthy expression
\[
     f^{}_{\II} (p^r) \; = \; \frac{(p+1)p^{r-4}}{(p-1)^3}\,
     \big( 4p^2 (2(p^2+1) + r (p^2-1)) + p^r (p^2+1)(r (p^4-1) 
     + p^4-4p^3-2p^2-4p+1)\big).
\]

Although this looks a bit nasty, the corresponding Dirichlet series 
has the nice form
\be
\begin{split}
  \varPhi^{}_{\II} (s) & \; = \; \varPhi(s) \varPhi(s-1)
     \; = \; \frac{\zeta^{}_{\II}(\frac{s}{2})
                   \zeta^{}_{\II}(\frac{s-1}{2})}
                  {\zeta^{}_L(2s) \zeta^{}_L(2s-2)}  \\
& \; =\; 1+\tfrac{25}{4^s}+\tfrac{36}{5^s}+\tfrac{100}{9^s}+\tfrac{288}{11^s}
 +\tfrac{440}{16^s}+\tfrac{800}{19^s}+\tfrac{900}{20^s}+\tfrac{960}{25^s}
 +\tfrac{1800}{29^s}+\tfrac{2048}{31^s}+\tfrac{2500}{36^s}+ \cdots
\end{split}
\ee which resembles the situation of the root lattice $D_4$ described
above, see \cite{BZ} for more. In particular, one can again use the
recursion \eqref{recursion} to calculate the coefficients $f^{}_{\II}
(m)$ directly from those of the icosahedral case. The last equation
also permits the determination of the asymptotic behaviour, compare
\cite[Appendix]{BM}. With the methods described in \cite[pp.\
29--31]{Wash}, applied to $L=\QQ(\tau)$, one can calculate the values
of $\zeta^{}_{L} (s)$ explicitly for $s\in \{2,4,6\}$.  This finally
gives that the number of CSMs of index $\le N$ is asymptotically
\[
     \frac{3^4 \, 5^7 \, 7}{268\, \pi^{12}}\sqrt{5}\log(\tau)\,
    \zeta^{}_{L} (3) N^3
     \; \simeq \; 0.19773\, N^3,
\]     
where $\zeta^{}_{L} (3) \simeq 1.02755$ has to be calculated numerically.

\section{Concluding remarks}

In this contribution, we have shown how the so-called coincidence
problem can be reformulated in a mathematical setting and then
solved algebraically in dimensions 2, 3 and 4. Various examples
have been treated explicitly, and it remains a simple exercise
to work out tables of all coincidence rotations with small
indices that could be relevant experimentally.

Rather obvious is the question for generalizations to higher
dimensions. One might hope that at least the root lattices could
be treated in full generality, but there are complications from
various sources. First of all, we do not have suitable generalizations
of quaternions at our disposal (and they proved extremely handy in
our treatment), and second, we depended on unique factorization
in one way or another -- and this does {\em not} generalize to
arbitrary lattices or modules.

Another obvious question emerges from the observation that we have
so far only dealt with {\em linear} isometries, while for various
reasons affine extensions are necessary, in particular for a
satisfactory formulation of the problem in the context of more
general Delone sets. Though some preliminary investigations exists
\cite{Pleasants}, more has to be done in this direction.

\appendix
\section{The number of sublattices of a given index}

Given a free Abelian group of rank $n$, one might like to know
how many different subgroups of (finite) index $m$ exist. Of course,
they are free and of rank $n$ again, but here we want to count
them separately, not up to isomorphism. Let us call that number
$f_n(m)$ and derive a recursion relation for it\footnote{This 
derivation is partially based on notes by P.~A.~B.\ Pleasants}. 
Since, for fixed $n$,
$f_n(m)$ is a multiplicative function in $m$, this will allow
the derivation of a closed formula both for $f_n(m)$ and
for its Dirichlet series generating function.

Since any free Abelian group of rank $n$ is isomorphic to $\ZZ^n$,
we can treat the latter case without loss of generality,
but with some benefit from the geometric setting.
Let $\gG$ be a sublattice of $\ZZ^n$ of index $m$.
Next, define a new lattice $\gL := \gG \cap \{x_n=0\}$
by intersection with an $(n\!-\!1)$-dimensional hyperplane.
Then, $\gL$ is a sublattice of $\ZZ^{n-1}$ of
finite index $[\ZZ^{n-1} : \gL] = d$, where we must
have $d \, | \, m$. 

At this point, we also know that $\gG$ can be generated
by $\gL$ and some vector $(\bs{y},m/d)$ with
$\bs{y}\in\ZZ^{n-1}$ (this is nothing but the completion
theorem for bases applied to this situation). Next, we observe
that we can actually calculate the number of lattices $\gG$
that give rise to the same $\gL$,
\[
   \big\lvert \{ \gG \mid [\ZZ^n:\gG]=m \;\mbox{ and }\;
        \gG\cap \{x_n=0\} = \gL \}\big\rvert
       \; = \; \big\lvert\{\mbox{choices for $\bs{y}$ (mod $\gL$)} \}
       \big\rvert   \; = \; d \, ,
\]
because $d$ is the number of residue classes of $\gL$
in $\ZZ^{n-1}$. Now, summation over all possibilities for
$\gL$ results in the following simple recursion formula
which is well-known, see \cite[$\S$ 63, A.~13 on p.~251]{Scheja}, 
but rather difficult to locate:
\be \label{mainrec}
   f_n(m) \; = \; \sum_{d\,|\,m} d \cdot f_{n-1}(d)
              \; = \; m \cdot \sum_{d\,|\,m} \frac{1}{m/d} 
              \cdot f_{n-1}(d) \, .
\ee

One can now derive a closed expression for $f_n(m)$, namely
\begin{prop} $\;\;$
      $ f_n(m) \; = \;  \sum_{d_1\cdot \ldots \cdot d_n = m}
           \,  d_1^0\cdot d_2^1\cdot \ldots\cdot d_n^{n-1} $ .
\end{prop}

Here, the sum runs over all $n$-tuples $(d_1,\ldots,d_n)$ of
positive integers subject to the restriction that
$d_1\cdot\ldots\cdot d_n = m$.

\noindent {\sc Proof}: It is clear that $f_1(m)\equiv1$.
From Eq.~(\ref{mainrec}), we get by induction 
\begin{eqnarray*}
   f_{n+1}(m) & = & \sum_{d\,|\,m}d\cdot f_n(d) \; = \; \sum_{d\,|\,m}
         \Big(d \sum_{d_1\cdot\ldots\cdot d_n=d} 
         d_1^0\cdot d_2^1\cdot\ldots\cdot d_n^{n-1} \Big) \\
   & = & \sum_{d\,|\,m} \hspace*{3mm} \sum_{d_1\cdot\ldots\cdot d_n=d}
         d_1^1\cdot d_2^2\cdot\ldots\cdot d_n^n   \\
   & = & \sum_{d_0\cdot d_1\cdot\ldots\cdot d_n=m}
         d_0^0\cdot d_1^1\cdot d_2^2\cdot\ldots\cdot d_n^n 
\end{eqnarray*}
which completes the argument, see \cite{Gruber,Zou} for
alternative approaches.  \qed

\smallskip
Let us determine a generating function for $f_n(m)$.
Due to the multiplica\-tivity of $f_n(m)$ in $m$, one would like
to have a Dirichlet series generating function. This can be
found as follows, using again the recursion relation \eqref{mainrec}.
\begin{eqnarray*}
   F_n(s+1)  & = & \sum_{m=1}^{\infty} \frac{f_n(m)}{m^{s+1}} 
         \; = \; \sum_{m=1}^{\infty} 
         \frac{\sum_{d\,|\,m} d \cdot f_{n-1}(d)}
              {m^{s+1}} \\
  & = & \sum_{m=1}^{\infty}
          \frac{\sum_{d\,|\,m} \frac{d}{m} \cdot f_{n-1}(d)}{m^s}
          \; = \; \sum_{m=1}^{\infty} \frac{1/m}{m^s} \cdot
          \sum_{\ell=1}^{\infty} \frac{f_{n-1}(\ell)}{\ell^s_{}} \\
  & = & \sum_{m=1}^{\infty} \frac{1}{m^{s+1}} \cdot F_{n-1}(s)
        \; = \; \zeta(s+1) \cdot F_{n-1}(s)\, .
\end{eqnarray*}
The middle line is the product formula for two Dirichlet series
generating functions, applied to our special case. From this
calculation, one gets the recursion 
\be 
     F_n(s) \; = \;  \zeta(s) \cdot F_{n-1}(s-1) \, .
\ee
Since there is only one sublattice of index $m$ for the case $n=1$,
we have $F_1(s) = \sum_{m=1}^{\infty} 1/m^s = \zeta(s)$ and
thus, by induction, one obtains
\begin{prop} $\;\;$
      $ F_n(s) \; = \; \zeta(s) \cdot \zeta(s-1) \cdot \ldots \cdot 
             \zeta(s-n+1)$ .     \qed
\end{prop}
In particular, this gives $F_2(s) = \zeta(s) \zeta(s-1)$,
which is the well-known generating function for the
divisor function $f_2^{}(m) = \sigma_1^{}(m) = \sum_{d\,|\,m}d$.

Let us close this appendix with a short remark on the asymptotic behaviour of
the coefficients. $F_n(s)$ has its rightmost pole at $s=n$, with residue
$r^{}_1=1$ (if $n=1$) and
$r^{}_n=\zeta(2)\cdot\zeta(3)\cdot\ldots\cdot\zeta(n)$ (if $n>1$). Then, the
number of sublattices with index $\le N$ is asymptotically given by
$r^{}_n\cdot N^n/n$, while the average number of sublattices with index $m$
grows like $m^{n-1}$.

\clearpage
\section*{Acknowledgements}

It is a pleasure to thank P.~A.~B.~Pleasants for his cooperation,
M.~Schlottmann for critically reading the manuscript, R.~V.~Moody,
J.~Roth, and A.~Weiss for helpful comments,
and the Fields Institute for Research in Mathematical Sciences
for financial support during a stay in fall 1995 where this manuscript
was written. I am also grateful to the MASCOS centre at the University of
Melbourne for support during a stay in 2006, where the manuscript
was revised and updated, and to P.~Zeiner for various suggestions
for improvements.

\end{document}